\documentclass[11pt]{article}

\usepackage{amsmath,amsfonts,amssymb,amsthm,latexsym}

\newcommand\qbi[3]{{{#1}\atopwithdelims[]{#2}}_{#3}}
\newcommand\bi[2]{{{#1}\atopwithdelims(){#2}}}

\newcommand{\la}{\lambda}
\def\C{\mathbb{ C}} \def\N{\mathbb{ N}}
 \def\Z{\mathbb{Z}}
\def\Q{\mathbb{Q}}

\newtheorem{theo}{Th\'eor\`eme}[section]
\newtheorem{prop}[theo]{Proposition}
\newtheorem{coro}[theo]{Corollaire}

\newtheorem{lem}[theo]{Lemme}
\newtheorem{rem}[theo]{Remarque}

\newenvironment{res}{\noindent{\bf R\'esum\'e. }}{\smallskip}
\newenvironment{abstr}{\noindent{\bf Abstract. }}{\smallskip}
\numberwithin{equation}{section}

\title{Irrationalit\'e aux entiers impairs positifs d'un $q$-analogue de la fonction z\^eta de Riemann}
\author{Fr\'ed\'eric Jouhet et Elie Mosaki}
\date{}

\begin{document}
\maketitle

\vspace{0.5 cm}

\noindent {\small \emph{2000 Mathematics Subject Classification}: Primary 11J72; Secondary 33D15.}

\vspace{0.7 cm}

\begin{res}
\small Dans cet article, nous nous int\'eressons \`a un $q$-analogue aux entiers positifs de la fonction z\^eta de Riemann, que l'on peut \'ecrire pour $s\in\N^*$ sous la forme $\zeta_q(s)=\sum_{k\geq 1}q^k\sum_{d|k}d^{s-1}$. Nous donnons une nouvelle minoration de la dimension de l'espace vectoriel sur $\Q$ engendr\'e, pour $1/q\in\Z\setminus\{-1;1\}$ et $A$ entier pair, par
$1,\zeta_q(3),\zeta_q(5),\dots,\zeta_q(A-1)$. Ceci am\'eliore un
r\'esultat r\'ecent de Krattenthaler, Rivoal et Zudilin
(\emph{S\'eries hyperg\'eom\'etriques basiques, $q$-analogues des
valeurs de la fonction z\^eta et s\'eries d'Eisenstein}, J. Inst.
Jussieu {\bf 5}.1 (2006), 53-79). En particulier notre r\'esultat a
pour cons\'equence le fait que pour $1/q\in\Z\setminus\{-1;1\}$, au moins l'un des nombres $\zeta_q(3),\zeta_q(5),\zeta_q(7),\zeta_q(9)$ est irrationnel.
\end{res}

\vspace{0.5 cm}

\begin{abstr}
\small In this paper, we focus on a
$q$-analogue of the Riemann z\^eta function at positive integers,
which can be written for $s\in\N^*$ by $\zeta_q(s)=\sum_{k\geq
1}q^k\sum_{d|k}d^{s-1}$. We give a new lower bound for the dimension of
the vector space over $\Q$ spanned, for
$1/q\in\Z\setminus\{-1;1\}$ and an even integer $A$, by
$1,\zeta_q(3),\zeta_q(5),\dots,\zeta_q(A-1)$. This improves a recent result of
Krattenthaler, Rivoal and Zudilin (\emph{S\'eries
hyperg\'eom\'etriques basiques, $q$-analogues des valeurs de la
fonction z\^eta et s\'eries d'Eisenstein}, J. Inst. Jussieu {\bf
5}.1 (2006), 53-79). In particular, a consequence of our result is that
for $1/q\in\Z\setminus\{-1;1\}$, at least one
of the numbers $\zeta_q(3),\zeta_q(5),\zeta_q(7),\zeta_q(9)$ is irrational.
\end{abstr}

\vspace{0.5 cm}

%%%%%%%%%%%%%%%%%%%%%%%%%%%%%%%%%%%%%%%%%%%%%%%%%%%%%%%%%%%%%%%%%%%%%%%%%%%%
\section{Introduction}
%%%%%%%%%%%%%%%%%%%%%%%%%%%%%%%%%%%%%%%%%%%%%%%%%%%%%%%%%%%%%%%%%%%%%%%%%%%%

L'\'etude de l'irrationalit\'e des valeurs de la fonction z\^eta
de Riemann $\zeta$ aux entiers impairs positifs est un probl\`eme classique en th\'eorie des nombres. Il est en
effet connu que l'expression des valeurs de $\zeta$ aux entiers
pairs positifs
$$\zeta(2m)=(-1)^{m-1}2^{2m-1}B_{2m}\frac{\pi^{2m}}{(2m)!}$$
permet d'affirmer, via la transcendance de $\pi$, due \`a
Lindemann, que chacun de ces nombres est transcendant (ici
$m\in\N^*$ et les nombres rationnels $B_m$ sont les nombres de
Bernoulli). En revanche, concernant l'\'etude aux entiers impairs
positifs, m\^eme si la transcendance est conjectur\'ee, le seul
r\'esultat significatif f\^ut pendant longtemps le th\'eor\`eme
d'Ap\'ery \cite{Ap} affirmant que \emph{$\zeta(3)$ est
irrationnel}. Puis r\'ecemment, Rivoal \cite{Ri1}, et Ball et
Rivoal \cite{BR} ont  eu l'id\'ee de consid\'erer les
valeurs de $\zeta$ aux entiers impairs positifs dans leur ensemble
plut\^ot qu'individuellement, ce qui leur permit de prouver qu'
\emph{il existe parmi les nombres $\zeta(2m+1)$, $m\in\N^*$, une
infinit\'e de nombres irrationnels}, en donnant la minoration pour
$A$ entier pair suffisamment grand :
$$\dim_{\Q}\left(\Q+\Q\zeta(3)+\dots+\Q\zeta(A-1)\right)\geq\frac{\log
A}{1+\log2}(1+\mbox{o}(1)).$$ 
La m\'ethode employ\'ee a conduit \`a des versions
quantitatives \cite{BR, KR, Ri1, Ri2}, jusqu'\`a l'article
r\'ecent de Zudilin \cite{Zu} dans lequel il est prouv\'e qu'
\emph{au moins l'un des nombres $\zeta(5), \zeta(7), \zeta(9), \zeta(11)$ est irrationnel}. Le lecteur int\'eress\'e pourra
aussi consulter le survol de
Fischler \cite{Fi} sur ce sujet.\\

Dans cet article, nous nous int\'eressons au $q$-analogue
normalis\'e de la fonction $\zeta$ consid\'er\'e d'abord dans
\cite{KKW} et \cite{Zu2}, puis plus r\'ecemment encore dans
\cite{KRZ}, et que l'on peut \'ecrire pour $s\in\N^*$ et $q$ un
nombre complexe tel que $|q|<1$~:
$$\zeta_q(s)=\sum_{k\geq
1}q^k\sum_{d|k}d^{s-1}=\sum_{k\geq1}k^{s-1}\frac{q^k}{1-q^k}\cdot$$ Le
terme de $q$-analogue est justifi\'e ici par la relation valide
pour $s\in\N^*\setminus\{1\}$ (voir par exemple \cite{KKW} ou
\cite{KRZ} pour une d\'emonstration) :
$$\lim_{q\to 1}(1-q)^s\zeta_q(s)=(s-1)!\zeta(s),$$
o\`u bien entendu $\zeta(s)=\sum_{k\geq 1}\frac{1}{k^s}$ est
l'expression pour $Re(s)>1$ de la fonction z\^eta de Riemann. L'un
des int\'er\^ets de ce $q$-analogue de $\zeta$ r\'eside dans le
fait que les valeurs de $\zeta_q$ aux entiers pairs positifs sont
reli\'ees aux formes modulaires et aux s\'eries d'Eisenstein
$E_{2m}(q)$ ($m\in \N^*$) \cite{Se} via la relation~:
$$E_{2m}(q)=1-\frac{4m}{B_{2m}}\zeta_q(2m).$$
Concernant la transcendance des valeurs de $\zeta_q$ aux entiers
pairs positifs, le r\'esultat d\'efinitif est cons\'equence de la structure de l'espace des formes modulaires sur $SL_2(\Z)$ \cite{Se} et d'un
th\'eor\`eme d'ind\'ependance alg\'ebrique sur les s\'eries
d'Eisenstein $E_2(q)$, $E_4(q)$ et $E_6(q)$ d\^u \`a Nesterenko \cite{Ne2}. En effet, on peut d\'eduire de cela que pour $m\in\N^*$ et $q$ alg\'ebrique (en
particulier $1/q\in\Z\setminus\{-1;1\}$), les
nombres $\zeta_q(2m)$ sont tous transcendants.

\medskip
Ceci conduit naturellement \`a se pencher sur le cas des valeurs
de $\zeta_q$ aux entiers impairs positifs. Remarquons tout d'abord
que malgr\'e l'analogie manifeste entre les r\'esultats de
transcendance des valeurs de $\zeta$ et $\zeta_q$ ($1/q\in
\Z\setminus \{-1;1\}$) aux entiers pairs positifs, il n'est
aujourd'hui possible d'affirmer l'irrationalit\'e de $\zeta_q(3)$
pour aucune valeur de $q$. En fait, seule l'irrationalit\'e de
$\zeta_q(1)$ est connue \cite{Bo} pour diverses valeurs de $q$. D'autre part, on sait depuis \cite{PV} que \emph{$1,\zeta_q(1),\zeta_q(2)$ sont
lin\'eairement ind\'ependants sur $\Q$ pour
$1/q\in\N\setminus\{1\}$}. Dans cette direction, le r\'esultat
principal de Krattenthaler, Rivoal et Zudilin dans \cite{KRZ}
affirme que pour $1/q\in\Z\setminus\{-1;1\}$ et $A$ entier pair~:
\begin{equation}\label{krzminor}
\dim_{\Q}\left(\Q+\Q\zeta_q(3)+\dots+\Q\zeta_q(A-1)\right)\geq
f(A),
\end{equation}
o\`u $$ f(A)=\max_{r\in\N\atop 1\leq r\leq
A/2}f(r;A)\;\;\;\;\mbox{avec}\;\;f(r;A):=\frac{4rA+A-4r^2}{\left(\frac{24}{\pi^2}+2\right)A+8r^2}\cdot$$
Cette minoration donne des informations asymptotiques via
l'\'equivalent 
$$f(A)\sim\frac{\pi}{2\sqrt{\pi^2+12}}\,\sqrt A\;\;\;\mbox{lorsque}\;\;A\to+\infty,$$
mais aussi quantitatives. En effet, il suffit de choisir une
valeur de $A\geq4$ la plus petite possible et donnant une
dimension sup\'erieure ou \'egale \`a 2 (l'id\'eal serait $A=4$,
ce qui montrerait l'irrationalit\'e de $\zeta_q(3)$). Cependant, il
s'av\`ere dans \cite{KRZ} que la valeur minimale exploitable est
$A=12$, ce qui fournit le r\'esultat suivant~: \emph{pour
$1/q\in\Z\setminus\{-1;1\}$, au moins l'un des nombres $\zeta_q(3),\zeta_q(5),\zeta_q(7),\zeta_q(9),\zeta_q(11)$ est irrationnel}.

\medskip
Le but de cet article est d'am\'eliorer \eqref{krzminor} et de raffiner
le r\'esultat quantitatif ci-dessus, en prouvant les deux
th\'eor\`emes suivants.
\begin{theo}
Pour $1/q\in\Z\setminus\{-1;1\}$ et tout entier pair $A\geq 4$, on a la
minoration~:
\begin{equation}\label{efminor}
\dim_{\Q}\left(\Q+\Q\zeta_q(3)+\dots+\Q\zeta_q(A-1)\right)\geq
g(A),
\end{equation}
o\`u
$$
g(A)=\max_{r\in\N\atop 1\leq r\leq
A/2}g(r;A)\;\;\;\;\mbox{avec}\;\;g(r;A):=\frac{4rA+A-4r^2}{\left(\frac{24}{\pi^2}+2\right)A-\frac{24}{\pi^2}+8r^2}\,,
$$
$g(A)$ v\'erifiant $\displaystyle \;g(A)\sim\frac{\pi}{2\sqrt{\pi^2+12}}\,\sqrt A\; $ lorsque $A\to+\infty$.
\end{theo}
\noindent On remarque ainsi qu'asymptotiquement, $g$ se comporte comme $f$ via
l'\'egalit\'e
$$\lim_{A\to+\infty}\frac{g(A)}{\sqrt A}=\lim_{A\to+\infty}\frac{f(A)}{\sqrt
A}=\frac{\pi}{2\sqrt{\pi^2+12}}\cdot$$ Cependant, pour toute valeur
fix\'ee de $A$, ce premier th\'eor\`eme am\'eliore la
minoration de \cite{KRZ} puisque $g(A)>f(A)$ (car $g(r;A)>f(r;A)$). Cette comparaison donne en particulier les in\'egalit\'es suivantes~:
\begin{eqnarray}
&&f(10)<1<g(10)=g(10;2)\simeq 1,001\label{dim1},\\
&&f(38)<g(38)<2<f(40)<g(40)\label{dim2},\\
&&f(86)<3<g(86)\label{dim3}.
\end{eqnarray}
La cons\'equence imm\'ediate de \eqref{dim1} est le Th\'eor\`eme~1.3 ci-dessous, qui est un raffinement de la version
quantitative de \cite{KRZ} d\'ej\`a mentionn\'ee. Les
in\'egalit\'es  \eqref{dim2} montrent que le Th\'eor\`eme~1.1
permet de retrouver, sans l'am\'eliorer, le r\'esultat suivant,
d\'ej\`a cons\'equence de \eqref{krzminor}~: \emph{pour
$1/q\in\Z\setminus\{-1;1\}$, il existe deux entiers impairs $j_1$
et $j_2$ tels que $3\leq j_1<j_2\leq 39$ et $1$, $\zeta_q(j_1)$ et
$\zeta_q(j_2)$ soient lin\'eairement ind\'ependants sur $\Q$}. En
revanche, \eqref{dim3} fournit une am\'elioration par rapport \`a
\eqref{krzminor}, qui peut s'\'ecrire~: 
\begin{coro}
Pour
$1/q\in\Z\setminus\{-1;1\}$, il existe trois entiers impairs
$j_1$, $j_2$ et $j_3$  tels que $3\leq j_1<j_2<j_3\leq
85$ et $1$, $\zeta_q(j_1)$, $\zeta_q(j_2)$ et $\zeta_q(j_3)$
soient lin\'eairement
ind\'ependants sur $\Q$.
\end{coro}
\medskip
\noindent Voici maintenant notre deuxi\`eme r\'esultat, qui est une
cons\'equence de \eqref{dim1}~:
\begin{theo}
Pour $1/q\in\Z\setminus\{-1;1\}$, au moins l'un des nombres $\zeta_q(3),\zeta_q(5)$, $\zeta_q(7),\zeta_q(9)$ est irrationnel.
\end{theo}
Il est int\'eressant de noter que la technique adopt\'ee dans
\cite{KRZ} est tout \`a fait parall\`ele (mais dans le monde des
$q$-analogues) \`a celle de \cite{Ri2}, o\`u il est d\'emontr\'e
qu'au moins l'un des nombres $\zeta(5),\zeta(7),\dots,\zeta(21)$ est irrationnel. Or les auteurs de \cite{KR} d\'emontrent la \emph{conjecture des d\'enominateurs} formul\'ee dans \cite{Ri3}, ce qui a pour cons\'equence le fait qu'au moins l'un des nombres
$\zeta(5),\zeta(7),\dots,\zeta(19)$ est irrationnel. C'est pourquoi le Th\'eor\`eme~1.3  n'est pas compl\`etement une surprise, les auteurs de \cite{KRZ} estimant \`a la fin de l'introduction qu'il est `probable' que l'on puisse prouver ce
r\'esultat, \`a condition de formuler et de d\'emontrer une certaine `$q$-conjecture des d\'enominateurs'. \\

Nous profitons de cette introduction pour donner les grandes
lignes de d\'emonstration du Th\'eor\`eme~1.1. Nous utilisons la proposition suivante, qui est un cas particulier du crit\`ere d'ind\'ependance lin\'eaire de Nesterenko \cite{Ne}~:
\begin{prop}
Soient un entier $N\geq 2$ et des r\'eels $v_1,\dots,v_N$.
Supposons qu'il existe $N$ suites d'entiers $(p_{j,n})_{n\geq0}$
et des r\'eels $\alpha_1$ et $\alpha_2$
avec $\alpha_2>0$ tels que~:\\
i)
$\displaystyle\lim_{n\to+\infty}\frac{1}{n^2}\log|p_{1,n}v_1+\dots+p_{N,n}v_N|=-\alpha_1$,\\
ii) pour tout $j\in\{1,\dots,N\}$, on a
$\displaystyle\limsup_{n\to+\infty}\frac{1}{n^2}\log|p_{j,n}|\leq\alpha_2$.\\
Alors la dimension du $\Q$-espace vectoriel engendr\'e par
$v_1,\dots,v_N$ v\'erifie~:
$$\dim_{\Q}\left(\Q v_1+\dots+\Q v_N\right)\geq 1+\frac{\alpha_1}{\alpha_2}\cdot$$
\end{prop}
\begin{rem}
\emph{Si en plus des hypoth\`eses de ce crit\`ere on connait un facteur
commun $\delta_n$ aux $p_{j,n}$, et si
$\displaystyle\lim_{n\to+\infty}\frac{1}{n^2}\log|\delta_n|$
existe et vaut $\delta$ ($<\alpha_2$), alors en consid\'erant les nouvelles suites
d'entiers $(p_{j,n}/\delta_n)_{n\geq0}$, on obtient~:
$$\dim_{\Q}\left(\Q v_1+\dots+\Q v_N\right)\geq
1+\frac{\alpha_1+\delta}{\alpha_2-\delta}\quad\left(\geq1+\frac{\alpha_1}{\alpha_2}\;\mbox{
si } \;\;\alpha_1>0\right).$$}
\end{rem}
\medskip
Afin d'exploiter le crit\`ere de Nesterenko dans notre contexte,
l'id\'ee consiste \`a analyser la s\'erie hyperg\'eom\'etrique
suivante (voir la partie~2 pour les notations)~:
\begin{equation*}
\tilde{S}_n(q):=(q)_n^{A-2r}\sum_{k\geq
1}(1-q^{2k+n})\frac{(q^{k-rn},q^{k+n+1})_{rn}}{(q^k)_{n+1}^A}q^{k(A-2r)n/2+kA/2-k},
\end{equation*}
o\`u $|q|\neq 1$, $A$ entier, $r\in\N^*$ et $A>2r$. Cette s\'erie
a \'et\'e sugg\'er\'ee, mais pas utilis\'ee, dans \cite{KRZ}, les
auteurs pr\'ef\'erant \'etudier une autre s\'erie, not\'ee
$S_n(q)$, pour prouver leurs r\'esultats. La premi\`ere \'etape
est une r\'e\'ecriture de $\tilde{S}_n(q)$, sous forme d'une
combinaison lin\'eaire en des $\zeta_q(2m+1)$, $m\in\N^*$ :
$$\tilde{S}_n(q)=\hat{P}_{0,n}(q)+\sum_{{j=3\atop j\,\mbox{\scriptsize
{impair}}}}^{A-1}\hat{P}_{j,n}(q)\zeta_q(j),$$ o\`u $|q|<1$, $A$ est
\emph{pair} et les $\hat{P}_{j,n}(q)$ sont \`a priori dans $\Q(q)$,
c'est-\`a dire des fractions rationnelles en $q$ (donc aussi en
$1/q$), \`a coefficients dans $\Q$. Dans un deuxi\`eme temps, on
cherche un d\'enominateur commun $D_n(q)$ \`a ces fractions
rationnelles en $1/q$, v\'erifiant~:
$$D_n(q)\hat{P}_{j,n}(q)\in\Z\left[\frac{1}{q}\right]\;\;\;\forall
j\in\{0,3,5,\dots,A-1\}.$$
Lorsque $1/q\in\Z\setminus\{-1;1\}$, la Proposition~1.4 appliqu\'ee \`a la combinaison lin\'eaire $D_n(q)\times \tilde{S}_n(q)$, ainsi que
les estimations asymptotiques de $\tilde{S}_n(q)$,
$\hat{P}_{j,n}(q)$ et $D_n(q)$, nous permettent de retrouver
\eqref{krzminor}.\\
 Notre am\'elioration, dont le r\'esultat est
donn\'e par le Th\'eor\`eme~1.1, se situe au niveau du
d\'enominateur commun $D_n(q)$ : nous formulons, puis
d\'emontrons, une  $q$-conjecture des d\'enominateurs (voir le Th\'eor\`eme~4.1) qui fournit un nouveau d\'enominateur commun $\tilde{D}_n(q)$ divisant $D_n(q)$. La minoration \eqref{efminor} du Th\'eor\`eme~1.1 est obtenue via l'estimation asymptotique de $\delta_n:=D_n(q)/\tilde{D}_n(q)$ (voir la Remarque~1.5 qui suit la Propri\'et\'e~1.4).
\begin{rem}
\emph{Pour d\'emontrer directement le Th\'eor\`eme~1.3, le crit\`ere de
Nestenrenko n'est pas n\'ecessaire. Il suffit en effet d'obtenir
une estimation asymptotique de la combinaison lin\'eaire \`a
coefficients entiers $\tilde {D}_n(q)\times \tilde{S}_n(q)$,
$1/q\in\Z\setminus\{-1;1\}$, en $\zeta_q(3)$, $\zeta_q(5)$,
$\zeta_q(7)$ et $\zeta_q(9)$ (avec les choix $A=10$ et $r=2$). Il
n'est donc pas n\'ecessaire de borner la hauteur des coefficients
de la combinaison lin\'eaire, ceci n'est utile que pour
l'ind\'ependance lin\'eaire.}
\end{rem}
\medskip
Cet article est organis\'e comme suit. La deuxi\`eme partie est destin\'ee \`a quelques notations concernant les $q$-s\'eries qui seront utiles ensuite. La troisi\`eme partie est
consacr\'ee \`a l'\'etude de la s\'erie $\tilde{S}_n(q)$ \'evoqu\'ee
ci-dessus. L'utilisation de $\tilde{S}_n(q)$ nous permet de
red\'emontrer la minoration \eqref{krzminor} de \cite{KRZ}, et de
d\'egager quelques lemmes cl\'es. Dans
cette m\^eme partie, nous expliquons par ailleurs comment le
nouveau d\'enominateur $\tilde{D}_n(q)$ permet d'obtenir le
Th\'eor\`eme~1.1. Dans la quatri\`eme et derni\`ere partie, nous
nous consacrons exclusivement \`a l'\'etude de $\tilde{D}_n(q)$ : nous
exprimons notre $q$-conjecture des d\'enominateurs (Th\'eor\`eme~4.1), et nous en donnons une d\'emonstration utilisant une
formule de transformation de s\'eries
hyperg\'eom\'etriques basiques due \`a Andrews \cite{An75, An84}.

%%%%%%%%%%%%%%%%%%%%%%%%%%%%%%%%%%%%%%%%%%%%%%%%%%%%%%%%%%%%%%%%%%%%%%%%%%%%
\section{Notations}
%%%%%%%%%%%%%%%%%%%%%%%%%%%%%%%%%%%%%%%%%%%%%%%%%%%%%%%%%%%%%%%%%%%%%%%%%%%%

Donnons comme annonc\'e ci-dessus quelques d\'efinitions et
notations issues du langage des $q$-s\'eries, que le lecteur
pourra retouver plus en d\'etails dans \cite{GR}. 

\medskip
Etant donn\'e 
un nombre complexe $q$ (la ``base'') tel que $|q|\neq 1$, on d\'efinit pour
tout r\'eel $a$ et tout entier $k\in\N$, le \emph{$q$-factoriel  montant} par~:
\begin{equation*}
(a)_k\equiv
(a;q)_k:=\left\{\begin{array}{l}1\;\;\mbox{si}\;\;k=0\\(1-a)\dots
(1-aq^{k-1})\;\;\mbox{si}\;\;k>0.
\end{array}\right.
\end{equation*}
La base $q$ peut \^etre omise lorsqu'il n'y a pas de confusion (en
notant $(a)_k$ pour $(a;q)_k$, etc), tout changement de base (par
exemple $q$ remplac\'e par $p=1/q$) sera pr\'ecis\'e
dans les paragraphes concern\'es. Pour des raisons pratiques, notons pour $k\in\N$~:
\begin{equation*}
(a_1,\ldots,a_m)_k:=(a_1)_k\times\cdots\times(a_m)_k.
\end{equation*}
\medskip
Rappelons aussi le \emph{coefficient $q$-binomial}~:
$$\qbi{n}{k}{q}=\qbi{n}{k,n-k}{q}:=\frac{(q)_n}{(q)_k(q)_{n-k}},$$
et plus g\'en\'eralement le \emph{coefficient $q$-multinomial}~:
$$\qbi{n}{k_1,\dots,k_l,n-k_1-\dots-k_l}{q}:=\frac{(q)_n}{(q)_{k_1}\dots(q)_{k_l}(q)_{n-k_1-\dots-k_l}},$$
qui sont des polyn\^omes en $q$, \`a coefficients
entiers (voir par exemple \cite{St}).

\medskip
\noindent Enfin, rappelons la notion
de \emph{s\'erie hyperg\'eom\'etrique basique ${}_{s+1}\phi_s$}~:
\begin{equation*}
{}_{s+1}\phi_s\!\left[\begin{matrix}a_0,a_1,\dots,a_s\\
b_1,\dots,b_{s}\end{matrix};q,z\right]:=
\sum_{k=0}^\infty\frac{(a_0,a_1,\dots,a_s)_k}{(q,b_1,\dots,b_{s})_k}z^k,
\end{equation*}
avec $a_j\in\C$ pour $0\leq j\leq s$, et $b_j\,q^k\neq 1$ pour
tout $k\in\N$ et $1\leq j\leq s$. La s\'erie converge toujours
pour $|z|<1$, et on dit que ${}_{s+1}\phi_s$ est~:
\begin{itemize}
\item \emph{bien \'equilibr\'ee} (well
poised) si $qa_0=a_1b_1=\dots=a_sb_s$ 
\item \emph{tr\`es bien
\'equilibr\'ee} (very well poised) si elle est bien \'equilibr\'ee
et de plus $a_1=q\sqrt{a_0}$.
\end{itemize}

%%%%%%%%%%%%%%%%%%%%%%%%%%%%%%%%%%%%%%%%%%%%%%%%%%%%%%%%%%%%%%%%%%%%%%%%%%%%
\section{Une s\'erie tr\`es bien \'equilibr\'ee}
%%%%%%%%%%%%%%%%%%%%%%%%%%%%%%%%%%%%%%%%%%%%%%%%%%%%%%%%%%%%%%%%%%%%%%%%%%%%

Reprenons la s\'erie d\'efinie dans l'introduction par~:
\begin{equation}\label{Sntilde}
\tilde{S}_n(q):=(q)_n^{A-2r}\sum_{k\geq
1}(1-q^{2k+n})\frac{(q^{k-rn},q^{k+n+1})_{rn}}{(q^k)_{n+1}^A}q^{k(A-2r)n/2+kA/2-k},
\end{equation}
pour $A$ entier pair, $r\in\N^*$ et $A-2r>0$. Remarquons
que cette s\'erie converge alors pour $|q|\neq1$. La s\'erie
$\tilde S_n(q)$ v\'erifie
\begin{equation}\label{tbe}
\tilde{S}_n(1/q)=-q^{n(r-1)}\tilde{S}_n(q),
\end{equation}
relation qui provient du choix de la puissance de $q$ dans le
sommande de \eqref{Sntilde}. La relation  (\ref{tbe}) va permettre
d'exprimer $\tilde{S}_n(q)$ comme combinaison lin\'eaire sur $\Q$ des
valeurs de $\zeta_q$ aux entiers impairs positifs seulement, alors
que l'on pourrait s'attendre \`a priori \`a voir appara\^itre
aussi les valeurs de $\zeta_q$ aux entiers pairs positifs. D'autre
part on peut \'ecrire
\begin{multline*}
\tilde{S}_n(q)=q^{(rn+1)((A-2r)n/2+A/2-1)}(1-q^{n+2rn+2})(q)_n^{A-2r}\frac{(q,q^{n+rn+2})_{rn}}{(q^{rn+1})_{n+1}^A}\\
\times{}_{A+4}\phi_{A+3}\!\left[\begin{matrix}a,q\sqrt a,-q\sqrt
a,q^{rn+1},\dots,q^{rn+1}\\
\sqrt a,-\sqrt
a,q^{(r+1)n+2},\dots,q^{(r+1)n+2}\end{matrix};q,q^{(A-2r)n/2+A/2-1}\right],
\end{multline*}
avec $a=q^{(2r+1)n+2}$, ce qui montre que $\tilde{S}_n(q)$ est
une s\'erie hyperg\'eom\'etrique basique tr\`es bien
\'equilibr\'ee.  Cette propri\'et\'e nous permettra de formuler
puis de d\'emontrer notre $q$-conjecture des d\'enominateurs dans
la partie~$4$ (voir le Th\'eor\`eme~4.1).

\medskip
\noindent Pour $|q|<1$, $A$ entier pair et $r\in\N^*$ tel que $A-2r>0$, nous allons successivement dans ce paragraphe d\'emontrer les~:
\begin{itemize}
\item  Lemme 3.2 (paragraphe 3.1) : on a
$$
\tilde{S}_n(q)=\hat{P}_{0,n}(q)+\sum_{{j=3\atop
j\,\mbox{\scriptsize
{impair}}}}^{A-1}\hat{P}_{j,n}(q)\zeta_q(j),$$ o\`u pour
$j\in\{0,3,5,\dots,A-1\}$, les $\hat{P}_{j,n}(q)$ sont des
fractions rationnelles  en $q$ qui seront explicit\'ees. 
\item
Lemme 3.5 (paragraphe 3.2) : si on pose
$d_n(q)=\mbox{ppcm}(q-1,\dots,q^n-1)$, alors pour
$j\in\{0,3,5,\dots,A-1\}$ et $\alpha=-A/8-r^2/2$, il existe des
r\'eels $\beta$ et $\gamma$ ne d\'ependant que de $A$ et $r$ tels
que
$$D_n(q)\hat{P}_{j,n}(q)\in\Z\left[\frac{1}{q}\right],\;\;\mbox{avec}\;\;D_n(q)=(A-1)!\,q^{\lfloor
\alpha n^2+\beta n+\gamma \rfloor}d_n(1/q)^{A},$$ o\`u
$\left\lfloor x\right\rfloor$ d\'esigne la partie enti\`ere de
$x$. \item Lemme 3.6 (paragraphe 3.3) : on a
$$\lim_{n\to+\infty}\frac{1}{n^2}\log|\tilde{S}_n(q)|=-\frac{1}{2}r(A-2r)\log|1/q|.$$
\item  Lemme 3.7 (paragraphe 3.3) : on a
$$\limsup_{n\to+\infty}\frac{1}{n^2}\log|\hat{P}_{j,n}(q)|\leq
\frac{1}{8}(A+4r^2)\log|1/q|,\;\; \forall
j\in\{0,3,5,\dots,A-1\}.$$ 
\item Lemme 3.8 (paragraphe 3.3) : on a
$$\lim_{n\to+\infty}\frac{1}{n^2}\log|D_n(q)|=\left(\frac{A}{8}+\frac{r^2}{2}+\frac{3A}{\pi^2}\right)\log|1/q|.$$
\end{itemize}

\medskip
\noindent 
Alors ces cinq lemmes permettent, pour $1/q\in \Z\setminus\{-1;1\}$,  d'appliquer la Proposition~1.4 \`a la combinaison lin\'eaire 
$$
D_n(q)\times\tilde{S}_n(q)=D_n(q)\hat{P}_{0,n}(q)+\sum_{{j=3\atop
j\,\mbox{\scriptsize
{impair}}}}^{A-1}D_n(q)\hat{P}_{j,n}(q)\zeta_q(j),$$ 
avec les valeurs
$$\alpha_1=-\left(\frac{A}{8}+\frac{r^2}{2}+\frac{3A}{\pi^2}-\frac{r}{2}(A-2r)\right)\log|1/q|$$
 et
$$\alpha_2=\left(\frac{A}{8}+\frac{r^2}{2}+\frac{3A}{\pi^2}+\frac{A}{8}+\frac{r^2}{2}\right)\log|1/q|=\left(\frac{A}{4}+r^2+\frac{3A}{\pi^2}\right)\log|1/q|,$$
ce qui implique
\begin{multline*}
\dim_\Q(\Q+\Q\zeta_q(3)+\Q\zeta_q(5)+\dots+\Q\zeta_q(A-1))\\
\geq1+\frac{\alpha_1}{\alpha_2}=\frac{4rA+A-4r^2}{\left(\frac{24}{\pi^2}+2\right)A+8r^2}\;,
\end{multline*}
red\'emontrant ainsi la minoration \eqref{krzminor}.

\medskip
Cependant, des calculs num\'eriques avec le logiciel Maple confirment (comme
le fait que $\tilde S_n(q)$ soit tr\`es bien \'equilibr\'ee le
laissait esp\'erer) que le d\'enominateur commun des
$\hat{P}_{j,n}(q)$ du Lemme~3.2 pourrait bien \^etre de la forme :
\begin{equation}\label{denominateurDn}
\tilde{D}_n(q)=(A-1)!\,q^{\lfloor \alpha n^2+\beta
n+\gamma\rfloor}d_n(1/q)^{A-1},\;\;\;\alpha=-\frac{A}{8}-\frac{r^2}{2},
\end{equation}
c'est-\`a-dire que l'on gagnerait une puissance de $d_n(1/q)$ par
rapport au choix $D_n(q)$. Ceci donnerait alors (voir la
Remarque~1.5 qui suit la proposition~1.4, avec $\delta=\lim_{n\to\infty}\frac{1}{n^2}\log d_n(1/q)$ et l'estimation
\eqref{dnasympt} de $d_n(1/q)$)~:
\begin{multline*}
\dim_\Q(\Q+\Q\zeta_q(3)+\Q\zeta_q(5)+\dots+\Q\zeta_q(A-1))\\
\geq1+\frac{\alpha_1+\delta}{\alpha_2-\delta}=\frac{4rA+A-4r^2}{\left(\frac{24}{\pi^2}+2\right)A-\frac{24}{\pi^2}+8r^2}\;,
\end{multline*}
ce qui d\'emontrerait le Th\'eor\`eme~1.1. Nous en d\'eduisons
donc que pour prouver le Th\'eor\`eme~1.1, il nous suffit de
montrer que le choix (\ref{denominateurDn}) est valide,
c'est-\`a-dire
\begin{equation}\label{Dntildeentier}
\tilde{D}_n(q)\hat{P}_{j,n}(q)\in\Z\left[\frac{1}{q}\right]\;\;\;\forall
j\in\{0,3,5,\dots,A-1\}\,,
\end{equation}
ce qui fera l'objet de la quatri\`eme partie.

%%%%%%%%%%%%%%%%%%%%%%%%%%%%%%%%%%%%%%%%%%%%%%%%%%%%%%%%%%%%
\subsection{Combinaisons lin\'eaires en les $\zeta_q(2j+1)$, $j\in\N^*$}

Posons
\begin{equation}\label{Rntilde}
\tilde{R}_n(T;q):=T^{(A-2r)n/2+A/2-2}q^{-An(n+1)/2}\frac{(q)_n^{A-2r}(q^{-rn}T,q^{n+1}T)_{rn}}{(T-1)^A
\dots(T-q^{-n})^A},
\end{equation}
de sorte que
$$\tilde{S}_n(q)=\sum_{k\geq 1}q^k(1-q^{2k+n})\tilde{R}_n(q^k;q).$$
Remarquons que le degr\'e en $T$ de la fraction rationnelle
$\tilde{R}_n(T;q)$ vaut $-(n+1)(A-2r)/2-r-2$ et est inf\'erieur ou
\'egal \`a $-3$ puisque $A>2r\geq2$. La d\'ecomposition de
cette fraction en \'el\'ements simples s'\'ecrit
$$
\tilde{R}_n(T;q)=\sum_{s=1}^A\sum_{j=0}^n\frac{\tilde{c}_{s,j,n}(q)}{(T-q^{-j})^s}=\sum_{s=1}^A\sum_{j=0}^n\frac{\tilde{d}_{s,j,n}(q)}{(1-Tq^{j})^s},
$$
avec $\tilde{d}_{s,j,n}(q)=(-1)^sq^{js}\tilde{c}_{s,j,n}(q)$ et
\begin{eqnarray}
\tilde{c}_{s,j,n}(q)&=&\frac{1}{(A-s)!}\left[\frac{d^{A-s}}{dT^{A-s}}\tilde{R}_n(T;q)(T-q^{-j})^A\right]_{T=q^{-j}}\label{ctilde}\\
&=&\frac{q^{-js}}{(A-s)!}\left[\frac{d^{A-s}}{du^{A-s}}\tilde{R}_n(uq^{-j};q)(u-1)^A\right]_{u=1}.\label{ctilde2}
\end{eqnarray}

\medskip
\noindent La d\'efinition (\ref{Rntilde}) conduit \`a la relation
$\tilde{R}_n(Tq^n;1/q)=q^{n(r-2)}\tilde{R}_n(T;q)$, qui a pour
cons\'equence pour tous $j\in\{0,\dots,n\}$ et $s\in\{1,\dots,A\}$
\begin{equation}\label{relationdj}
\tilde{d}_{s,n-j,n}(1/q)=q^{n(r-2)}\tilde{d}_{s,j,n}(q),
\end{equation}
ou de fa\c con \'equivalente
\begin{equation}\label{relationcj}
\tilde{c}_{s,n-j,n}(1/q)=q^{n(s+r-2)}\tilde{c}_{s,j,n}(q).
\end{equation}
\begin{rem}
\emph{La relation \eqref{relationdj} est un peu diff\'erente de celle
prouv\'ee dans \cite{KRZ}, sauf dans le cas $r=1$, ce qui
n'est pas une surprise car la s\'erie $\tilde{S}_n(q)$ coincide pour
$r=1$ avec la s\'erie $S_n^{[1]}(q):=S_n(q)-S_n(1/q)$ utilis\'ee
dans \cite{KRZ}.}
\end{rem}

\medskip

Nous aurons besoin dans ce qui suit des \emph{nombres de Stirling
de premi\`ere esp\`ece sans signe} (voir \cite{St}), qui sont des
nombres entiers not\'es $c(s,j)$ (o\`u $s$ et $j$ sont deux
entiers tels que $1\leq j\leq s$) et d\'efinis par
$$x(x+1)\dots(x+s-1)=\sum_{j=1}^sc(s,j)x^j.$$
\begin{lem}
 On a pour $|q|<1$, $A$ pair et $r\in\N^*$ tel que $A-2r>0$~:
\begin{equation}\label{formelineaire}
\tilde{S}_n(q)=\hat{P}_{0,n}(q)+\sum_{{j=3\atop
j\,\mbox{\scriptsize {impair}}}}^{A-1}\hat{P}_{j,n}(q)\zeta_q(j),
\end{equation}
o\`u pour $j=3,5,\dots,A-1$,
\begin{eqnarray}
\hat{P}_{0,n}(q)&:=&\tilde{P}_{0,n}(1,q)-q^{-n(r-1)}\tilde{P}_{0,n}(1,1/q)-
\left[\frac{d}{dz}\tilde{P}_{1,n}(z,q)\right]_{z=1}\!\!\label{P0chap},\\
\hat{P}_{j,n}(q)&:=&\sum_{s=j}^A\frac{2c(s-1,j-1)}{(s-1)!}\tilde{P}_{s,n}(1,q)\label{Pjchap},
\end{eqnarray}
et
\begin{eqnarray}
\tilde{P}_{0,n}(z,q)&:=&\sum_{s=1}^A\sum_{j=1}^n\sum_{k=1}^j(-1)^s\frac{q^{k-j(1-s)}}{(1-q^k)^s}\,\tilde{c}_{s,j,n}(q)z^{j-k}\label{P0tilde},\\
\tilde{P}_{s,n}(z,q)&:=&(-1)^s\sum_{j=0}^nq^{j(s-1)}\tilde{c}_{s,j,n}(q)z^j\label{Pjtilde}.
\end{eqnarray}
\end{lem}
\begin{proof}[D\'emonstration]
Reprenons les fonctions interm\'ediaires de \cite{KRZ} d\'efinies
par $Z_s(z;q):=\sum_{k\geq1}\frac{q^k}{(1-q^k)^s}z^{-k}$ et
$Z_s(q):=Z_s(1;q)$, v\'erifiant notamment pour $s\geq 2$
\begin{equation}\label{zs-stirling}
Z_s(q)-Z_s(1/q)=\frac{2}{(s-1)!}\sum_{{j=3\atop
j\,\mbox{\scriptsize {impair}}}}^{s}c(s-1,j-1)\zeta_q(j),
\end{equation}
et convergeant pour tout $s\geq 1$ d\`es que $|q|<|z|$. Nous
utiliserons le fait que $Z_s(1/z;1/q)$ converge pour tout
$|z|<|q|^{1-s}$. Posons alors pour $|q|<|z|$
$$
{\cal{S}}_n(z;q):=\sum_{k\geq 1}q^k\tilde{R}_n(q^k;q)z^{-k}.
$$
Il n'est pas difficile de voir que ${\cal{S}}_n(z;q)$ converge
d\`es que $|z|>|q|^{(A-2r)n/2+A/2-1}$ et ${\cal{S}}_n(1/z;1/q)$
converge d\`es que $|z|<|q|^{-(A-2r)n/2-A/2-1}$. Les deux s\'eries ${\cal{S}}_n(z;q)$ et ${\cal{S}}_n(1/z;1/q)$ convergent donc simultan\'ement pour $|q|<|z|\leq1$ lorsque $A\geq4$. Un calcul simple
montre que
\begin{equation}\label{**}
{\cal{S}}_n(1;q)-q^{-n(r-1)}{\cal{S}}_n(1;1/q)=\tilde{S}_n(q).
\end{equation}

\medskip
\noindent 
Par ailleurs, en utilisant la d\'ecomposition en \'el\'ements
simples de $\tilde{R}_n(T;q)$, on obtient
\begin{equation}\label{devcalSn}
{\cal{S}}_n(z;q)=\tilde{P}_{0,n}(z,q)+\sum_{s=1}^A\tilde{P}_{s,n}(z,q)Z_s(z;q),
\end{equation}
o\`u $\tilde P_{0,n}(z,q)$ et $\tilde P_{s,n}(z,q)$ sont des polyn\^omes en
$z$ d\'efinis par \eqref{P0tilde} et \eqref{Pjtilde}. Or, en
reprenant la d\'efinition (\ref{Pjtilde}), on s'aper\c coit que la
relation (\ref{relationdj}) se traduit pour $s\geq1$ par
$\tilde{P}_{s,n}(1/z,1/q)=z^{-n}q^{n(r-1)}\tilde{P}_{s,n}(z,q)$,
ce qui nous donne l'id\'ee d'\'etudier maintenant la s\'erie
\begin{equation}\label{defcalSntilde}
\tilde{{\cal{S}}}_n(z;q):={\cal{S}}_n(z;q)-z^{n}q^{-n(r-1)}{\cal{S}}_n(1/z;1/q),
\end{equation}
avec la condition $|q|<|z|<1$ assurant la convergence. Le
d\'eveloppement (\ref{devcalSn}), puis la relation ci-dessus entre
$\tilde{P}_{s,n}(1/z,1/q)$ et $\tilde{P}_{s,n}(z,q)$, permettent
alors d'\'ecrire pour $|q|<|z|<1$
\begin{multline}\label{devcalSntilde}
\tilde{{\cal{S}}}_n(z;q)=\tilde{P}_{0,n}(z,q)-q^{-n(r-1)}z^{-n}\tilde{P}_{0,n}(1/z,1/q)\\
+\sum_{s=1}^A\tilde{P}_{s,n}(z,q)(Z_s(z;q)-Z_s(1/z;1/q)).
\end{multline}
Il ne reste plus qu'\`a faire tendre $z$ vers 1 dans
(\ref{devcalSntilde}) pour obtenir le Lemme~3.2, via la
d\'efinition (\ref{defcalSntilde}) et les \'egalit\'es
(\ref{zs-stirling}) et (\ref{**}), \'etant entendu (comme dans
\cite{KRZ}) que
$$\lim_{z\to1}\tilde{P}_{1,n}(z,q)(Z_1(z;q)-Z_1(1/z;1/q))=-\left[\frac{d}{dz}\tilde{P}_{1,n}(z,q)\right]_{z=1}.$$
\end{proof}

%%%%%%%%%%%%%%%%%%%%%%%%%%%%%%%%%%%%%%%%%%%%%%%%%%%%%%%%%%%%
\subsection{Propri\'et\'es arithm\'etiques des $\hat{P}_{j,n}$}

En vue de trouver un d\'enominateur commun $D_n(q)$ aux
coefficients $\hat{P}_{j,n}(q)\in\Q(q)$, rappelons que
$d_n(q)\in\Z[q]$ est le polyn\^ome unitaire, de plus petit degr\'e
et multiple commun de $1-q,1-q^2,\dots,1-q^n$. On d\'emontre alors
les deux lemmes suivants, qui permettront de d\'eduire le Lemme
$3.5$.
\begin{lem}
Pour tous $s\in\{1,\dots,A\}$ et $j\in\{0,\dots,n\}$, on a :
$$d_n\left(1/q\right)^{A-s}\tilde{c}_{s,j,n}(q)\in\mathbb{Z}\left[q;\frac{1}{q}\right].$$
\end{lem}
\begin{proof}[D\'emonstration]
\'Ecrivons (\ref{ctilde}) sous la forme
\begin{equation}\label{ctilderecrit}
\tilde{c}_{s,j,n}(q)=\frac{1}{(A-s)!}\left[\frac{d^{A-s}}{dT^{A-s}}\tilde{V}_n(T;q)\right]_{T=q^{-j}},
\end{equation}
avec
\begin{eqnarray}
\tilde{V}_n(T;q)&:=&\tilde{R}_n(T;q)(T-q^{-j})^A\nonumber\\
&=&(q)_n^{A-2r}T^{(A-2r)n/2+A/2-2}q^{-An(n+1)/2}\nonumber\\
&&\hskip
2cm\times\frac{(q^{-rn}T,q^{n+1}T)_{rn}}{(T-1)^A\dots(T-q^{-n})^A}(T-q^{-j})^A.\label{V}
\end{eqnarray}
On regroupe les termes de $\tilde{V}_n(T;q)$ de la mani\`ere
suivante :
$$\tilde{V}_n(T;q)=q^{an^2+bn+c}T^{A/2-2}F(T)^{A/2-r}G(T)^{A/2-r}\prod_{l=1}^rH_l(T)I_l(T),
$$
o\`u $a$, $b$ et $c$ sont des entiers d\'ependant uniquement de
$A$ et $r$, les fonctions $F$, $G$, $H_l$ et $I_l$ \'etant celles
d\'efinies dans \cite{KRZ}. Donc en utilisant leurs
d\'ecompositions en \'el\'ements simples (voir \cite{KRZ}), on
s'aper\c coit que si l'on note $U$ n'importe laquelle de ces
fonctions, ou m\^eme $T\mapsto T^{A/2-2}$, alors
$$\frac{d_n\left(1/q\right)^{\mu}}{\mu!}\left[\frac{d^{\mu}}{dT^{\mu}}U(T)\right]_{T=q^{-j}}\in\mathbb{Z}\left[q;\frac{1}{q}\right]\;\;\;\forall\mu\in\N.$$
On conclut en utilisant (\ref{ctilderecrit}) et en appliquant la
formule de Leibniz de d\'erivation $(\mu=A-s)$-i\`eme d'un produit
de fonctions.
\end{proof}

\medskip
\noindent 

\begin{lem}
Soit $\alpha=-A/8-r^2/2$. Il existe alors $\beta'$ et $\gamma'$
r\'eels d\'ependant uniquement de $A$ et $r$ tels que pour tous $(s,j)\in\{1,\dots,A\}\times\{0,\dots,n\}$~:
$$\lim_{q\to+\infty}q^{\alpha n^2+\beta'
n+\gamma'}\tilde{c}_{s,j,n}(q)<\infty.$$
\end{lem}
\begin{proof}[D\'emonstration]
 Reprenons l'expression (\ref{ctilderecrit}) de la d\'emonstration pr\'ec\'edente,
et posons
$\tilde{v}_n(T;q):=\frac{d}{dT}\tilde{V}_n(T;q)/\tilde{V}_n(T;q)$
la d\'eriv\'ee logarithmique en $T$ de $\tilde{V}_n(T;q)$. Comme
dans \cite{KRZ}, la formule de d\'erivation de Fa\`a di Bruno
donne alors pour tout $\mu\in\N$
\begin{equation}\label{FaaVn}
\frac{1}{\mu!}\frac{d^{\mu}}{dT^{\mu}}\tilde{V}_n(T;q)=\sum_{k_1+\dots+\mu
k_\mu=\mu}\frac{\tilde{V}_n(T;q)}{k_1!\dots
k_\mu!}\prod_{l=1}^\mu\left(\frac{1}{l!}\frac{d^{l-1}}{dT^{l-1}}\tilde{v}_n(T;q)\right)^{k_l}.
\end{equation}
Or par d\'efinition
\begin{eqnarray*}
\tilde{v}_n(T;q)&=&\frac{d}{dT}(\log\tilde{V}_n(T;q))\\
&=&\frac{(A-2r)n/2+A/2-2}{T}+\sum_{i=1}^{rn}\frac{1}{T-q^i}\\
&&\hskip 3cm
+\sum_{i=n+1}^{rn+n}\frac{1}{T-q^{-i}}-A\sum_{{i=0\atop i\neq
j}}^{n}\frac{1}{T-q^{-i}},
\end{eqnarray*}
donc pour $l\in\N$
\begin{multline*}
\frac{(-1)^{l-1}}{(l-1)!}\frac{d^{l-1}}{dT^{l-1}}\tilde{v}_n(T;q)=\frac{(A-2r)n/2+A/2-2}{T^l}+\sum_{i=1}^{rn}\frac{1}{(T-q^i)^l}\\
+\sum_{i=n+1}^{rn+n}\frac{1}{(T-q^{-i})^l}-A\sum_{{i=0\atop i\neq
j}}^{n}\frac{1}{(T-q^{-i})^l}\cdot
\end{multline*}
Or on peut \'ecrire ceci sous la forme
\begin{multline*}
\frac{(-1)^{l-1}}{(l-1)!}\frac{d^{l-1}}{dT^{l-1}}\tilde{v}_n(T;q)=\frac{(A-2r)n/2+A/2-2}{T^l}+\sum_{i=1}^{rn}\left(\frac{q^{-i}}{Tq^{-i}-1}\right)^l\\
+\sum_{i=n+1}^{rn+n}\frac{1}{(T-q^{-i})^l}-A\sum_{i=0}^{j-1}\left(\frac{q^{i}}{Tq^{i}-1}\right)^l-A\sum_{i=j+1}^{n}\frac{1}{(T-q^{-i})^l},
\end{multline*}
ce qui implique que $\forall j\in\{0,\dots,n\}$,
$\displaystyle\lim_{q\to+\infty}q^{-jl}\left[\frac{d^{l-1}}{dT^{l-1}}\tilde{v}_n(T;q)\right]_{T=q^{-j}}<\infty$.
On en d\'eduit que pour $k_1+\dots+\mu k_\mu=\mu$,
\begin{equation}\label{E}
\lim_{q\to+\infty}q^{-j\mu}\left[\prod_{l=1}^\mu\left(\frac{1}{l!}\frac{d^{l-1}}{dT^{l-1}}\tilde{v}_n(T;q)\right)^{k_l}\right]_{T=q^{-j}}<\infty.
\end{equation}
Par ailleurs la puissance dominante de $q$ apparaissant dans $\tilde{V}_n(q^{-j};q)$ d\'efini gr\^ace \`a \eqref{V}
est de la forme~:
$$j(nA/2-A+2)-j^2A/2+rn(rn-1)/2\,.$$

\medskip
\noindent 
Il suffit maintenant de choisir $\mu=A-s$ dans (\ref{FaaVn}), puis \`a l'aide de (\ref{ctilderecrit}) et de \eqref{E} on obtient~:
$$\lim_{q\to+\infty}q^{-j(A-s)-\left(j(nA/2-A+2)-j^2A/2+rn(rn-1)/2\right)}\tilde{c}_{s,j,n}(q)<\infty\;.$$
On conclut alors ais\'ement puisque pour tous
$(s,j)\in\{1,\dots,A\}\times\{0,\dots,n\}$,
$$-j(A-s)-\left(j(nA/2-A+2)-j^2A/2+rn(rn-1)/2\right)\geq \alpha n^2+\beta' n+\gamma'\;,$$
avec $\alpha=-A/8-r^2/2$, $\beta'=(r-1)/2$, $\gamma'=-1/(2A)$ (cette
borne inf\'erieure est obtenue pour $s=1$ et $j=n/2+1/A$).
\end{proof}

\medskip
 
Le Lemme~3.4 implique que $\lim_{q\to+\infty}q^{\alpha n^2+\beta'
n+\gamma'}d_n(1/q)^{A-s}\tilde{c}_{s,j,n}(q)<\infty$ pour tous
$(s,j)\in\{1,\dots,A\}\times\{0,\dots,n\}$ puisque
$\lim_{q\to+\infty}d_n(1/q)=d_n(0)=\pm1$. En utilisant le
Lemme~3.3 on obtient donc que pour $\alpha=-A/8-r^2/2$ il existe
$\beta'$ et $\gamma'$ r\'eels ne d\'ependant que de $A$ et $r$, tels
que pour tous $ (s,j)\in\{1,\dots,A\}\times\{0,\dots,n\}$
\begin{equation}\label{100}
q^{\lfloor\alpha n^2+\beta'
n+\gamma'\rfloor}d_n\left(1/q\right)^{A-s}\tilde{c}_{s,j,n}(q)\in\Z\left[\frac{1}{q}\right].
\end{equation}
D'apr\`es les expressions (\ref{P0chap})-(\ref{Pjtilde}), l'\'equation \eqref{relationcj}, la d\'efinition de $d_n(1/q)$ et le fait que
\begin{equation}\label{*}
\left[\frac{d}{dz}\tilde{P}_{1,n}(z,q)\right]_{z=1}=-\sum_{j=0}^nj\,\tilde{c}_{1,j,n}(q),
\end{equation}
on d\'eduit de \eqref{100} le lemme suivant (de simples calculs montrent que les valeurs $\beta=\beta'-A+1$ et $\gamma=\gamma'+A-2$ conviennent)~:
\begin{lem}
Pour $\alpha=-A/8-r^2/2$, il existe $\beta$ et $\gamma$ r\'eels ne
d\'ependant que de $A$ et $r$ tels que~:
\begin{equation}\label{denominatpjchap}
(A-1)!\,q^{\lfloor\alpha n^2+\beta
n+\gamma\rfloor}d_n(1/q)^{A-j}\hat{P}_{j,n}(q)\in\Z\left[\frac{1}{q}\right]\;\;\forall
j\in\{3,5,\dots,A-1\}
\end{equation}
et
\begin{equation}\label{denominatp0chap}
q^{\lfloor\alpha n^2+\beta
n+\gamma\rfloor}d_n(1/q)^A\hat{P}_{0,n}(q)\in\Z\left[\frac{1}{q}\right].
\end{equation}
Ainsi, en posant
\begin{equation}\label{DenominateurcommunDn}
D_n(q)=(A-1)!\,q^{\lfloor\alpha n^2+\beta
n+\gamma\rfloor}d_n(1/q)^{A},
\end{equation}
on obtient~:
$$D_n(q)\hat{P}_{j,n}(q)\in\Z\left[\frac{1}{q}\right]\;\;\forall
j\in\{0,3,5,\dots,A-1\}.$$
\end{lem}

%%%%%%%%%%%%%%%%%%%%%%%%%%%%%%%%%%%%%%%%%%%%%%%%%%%%%%%%%%%%
\subsection{Estimations asymptotiques}

On \'evalue maintenant asymptotiquement $\tilde{S}_n(q)$, puis les
coefficients
$\hat{P}_{j,n}(q)$ de (\ref{formelineaire}), et enfin $D_n(q)$. Fixons $A$ entier pair et $r\in\N^*$ tel que $A-2r>0$.\\

Commen\c cons par l'estimation asymptotique de $\tilde{S}_n(q)$,
donn\'ee par le lemme suivant~:
\begin{lem}
Pour tout $|q|<1$, on a :
$$\lim_{n\to+\infty}\frac{1}{n^2}\log|\tilde{S}_n(q)|=-\frac{1}{2}r(A-2r)\log|1/q|.$$
\end{lem}
\begin{proof}[D\'emonstration]
On note $\rho_k(q)=q^k(1-q^{2k+n})\tilde{R}_n(q^k;q)$, de sorte
que $\tilde{S}_n(q)=\sum_{k\geq 1}\rho_k(q)$. Par la d\'efinition \eqref{Rntilde}
 de $\tilde{R}_n(T;q)$, il apparait clairement que
$\rho_k(q)=0\Leftrightarrow k\in\{0,\dots,rn\}$. Or on a pour
$k\geq rn+1$
$$\frac{\rho_{k+1}(q)}{\rho_k(q)}=q^{(A-2r)n/2+A/2-1}\frac{1-q^{2k+n+2}}{1-q^{2k+n}}\frac{1-q^{1+k+n+rn}}{1-q^{k-rn}}\left(\frac{1-q^{k}}{1-q^{k+n+1}}\right)^{A+1}.$$
Comme $A-2r>0$, $|q|<1$ et $k\geq rn+1$, on a donc pour $n$ grand la
majoration uniforme en $k$~:
$$\left|\frac{\rho_{k+1}(q)}{\rho_k(q)}\right|\leq|q|^{(A-2r)n/2}\left(\frac{1+|q|}{1-|q|}\right)^{A+3}<\frac{1}{3},$$
ce qui permet comme dans \cite{KRZ} d'\'ecrire l'encadrement
$$\frac{1}{2}|\rho_{rn+1}(q)|\leq|\tilde{S}_n(q)|\leq\frac{3}{2}|\rho_{rn+1}(q)|.$$
Or
$$\rho_{rn+1}(q)=(1-q^{2rn+n+2})(q)_n^{A-2r}\frac{(q,q^{(r+1)n+2})_{rn}}{(q^{rn+1})_{n+1}^A}q^{(rn+1)((A-2r)n/2+A/2-1)},$$
donc on obtient
$$\lim_{n\to+\infty}\frac{1}{n^2}\log|\tilde{S}_n(q)|=\lim_{n\to+\infty}\frac{1}{n^2}\log|\rho_{rn+1}(q)|=-\frac{1}{2}r(A-2r)\log|1/q|.$$
\end{proof}

Donnons maintenant l'estimation asymptotique des coefficients
$\hat{P}_{j,n}(q)$ de (\ref{formelineaire}), \`a l'aide du lemme
suivant~:
\begin{lem}
Pour tout $j\in\{0,3,5,\dots,A-1\}$ et $|q|<1$, on a :
$$\limsup_{n\to+\infty}\frac{1}{n^2}\log|\hat{P}_{j,n}(q)|\leq
\frac{1}{8}(A+4r^2)\log|1/q|.$$
\end{lem}
\begin{proof}[D\'emonstration]
Remarquons tout d'abord que pour des nombres complexes $a_{i,n}$
($0\leq i\leq n$),
$$\left(\forall
i\in\{0,\dots,n\},\,\limsup_{n\to+\infty}\frac{1}{n^2}\log|a_{i,n}|\leq
c\right)\,\Rightarrow\limsup_{n\to+\infty}\frac{1}{n^2}\log\left|\sum_{i=0}^na_{i,n}\right|\leq
c.$$
 Ceci montre, via les d\'efinitions des $\hat{P}_{j,n}(q)$ donn\'ees par
(\ref{P0chap})-(\ref{Pjtilde}), qu'il suffit de prouver que
l'estimation du lemme est valide pour les coefficients
$\tilde{d}_{s,j,n}(q)=(-1)^sq^{js}\tilde{c}_{s,j,n}(q)$,
uniform\'ement en $j$ et $s$. Fixons maintenant $j\in\{0,\dots,n\}$ et
$\eta=(1-|q|)/2>0$; la formule de Cauchy appliqu\'ee \`a
(\ref{ctilde}) donne alors
$$\tilde{d}_{s,j,n}(q)=-\frac{1}{2i\pi}\int_{\cal
C}\tilde{R}_n(Tq^{-j};q)(1-T)^{s-1}dT,$$ o\`u ${\cal C}$ d\'esigne
le cercle de centre 1 et de rayon $\eta$. Reprenons donc
l'expression (\ref{Rntilde}), qui conduit \`a
\begin{multline*}
\tilde{R}_n(Tq^{-j};q)(1-T)^{s-1}=q^{-j((A-2r)n/2+A/2-2)}T^{(A-2r)n/2+A/2-2}\\
\times(q)_n^{A-2r}(1-T)^{s-1}\frac{(q^{-rn-j}T,q^{n-j+1}T)_{rn}}{(Tq^{-j})_{n+1}^A}\cdot
\end{multline*}
Ceci peut s'\'ecrire apr\`es quelques transformations
\'el\'ementaires
\begin{multline*}
\tilde{R}_n(Tq^{-j};q)(1-T)^{s-1}=q^{Aj^2/2-Anj/2+2j-rn(rn+1)/2}T^{A(n-2j)/2+A/2-2}\\
\times(-1)^{Aj+rn}(q)_n^{A-2r}(1-T)^{s-A-1}\frac{(q^{j+1}/T,q^{n-j+1}T)_{rn}}{(q/T)_{j}^A(qT)_{n-j}^A}\cdot
\end{multline*}
En vue de majorer cette expression pour $T\in{\cal C}$, on reprend
maintenant les encadrements de \cite{KRZ} valables pour
$(a,b)\in\N^*\times\N$,  $T\in{\cal C}$ et $\eta=(1-|q|)/2$ :
$$0<(|q|(1+\eta);|q|)_\infty\leq|(q^aT)_b|\leq(-(1+\eta);|q|)_\infty$$
et
$$0<(|q|/(1-\eta);|q|)_\infty\leq|(q^a/T)_b|\leq(-1/(1-\eta);|q|)_\infty,$$
puis
$$|T^{A(n-2j)/2+A/2-2}|\leq (\max(1+\eta;1/(1-\eta))^{An/2}(1+\eta)^{A/2-2}$$
et
$$|(q)_n|\leq(-|q|;|q|)_\infty\;\;\mbox{et}\;\;|1-T|^{s-A-1}\leq1/\eta^{A+1}.$$

\medskip
\noindent 
Il ne reste donc plus qu'\`a majorer la puissance de $q$ dans
l'expression de $\tilde{R}_n(Tq^{-j};q)(1-T)^{s-1}$ ci-dessus. Ceci se fait en remarquant que la fonction $j\mapsto
Aj^2/2-Anj/2+2j-rn(rn+1)/2$ atteint son maximum en $j=n/2-2/A$, et
cette valeur maximale est de la forme $-An^2/8-r^2n^2/2+\la
n+\mu$, $\lambda$ et $\mu$ \'etant des r\'eels ne d\'ependant que
de $A$ et $r$. Tout cela conduit \`a la majoration
$$|\tilde{d}_{s,j,n}(q)|\leq\tilde{c}_0\times|q|^{-(A+4r^2)n^2/8},$$
o\`u $\tilde{c}_0$ ne d\'epend ni de $j$ ni de $s$, et v\'erifie
$\displaystyle\lim_{n\to+\infty}\tilde{c}_0^{\;1/n^2}=1$, ce qui
permet de conclure.
\end{proof}
Finalement, donnons l'estimation asymptotique de $D_n(q)$ d\'efini
en (\ref{DenominateurcommunDn})~:
\begin{lem}
Pour tout $|q|<1$ on a
$$\lim_{n\to+\infty}\frac{1}{n^2}\log|D_n(q)|=\left(\frac{A}{8}+\frac{r^2}{2}+\frac{3A}{\pi^2}\right)\log|1/q|.$$
\end{lem}
\begin{proof}[D\'emonstration]
Pour $|q|<1$ on a l'estimation (voir \cite{BV} et \cite{VA})
\begin{equation}\label{dnasympt}
\lim_{n\to+\infty}\frac{1}{n^2}\log|d_n(1/q)|=\frac{3}{\pi^2}\log|1/q|,
\end{equation} donc la conclusion est imm\'ediate \`a l'aide des expressions de
$D_n(q)$ et de $\alpha$ donn\'ees par
(\ref{DenominateurcommunDn}).
\end{proof}

%%%%%%%%%%%%%%%%%%%%%%%%%%%%%%%%%%%%%%%%%%%%%%%%%%%%%%%%%%%%%%%%%%%%%%%%%%%%
\section{D\'emonstration du Th\'eor\`eme~1.1}
%%%%%%%%%%%%%%%%%%%%%%%%%%%%%%%%%%%%%%%%%%%%%%%%%%%%%%%%%%%%%%%%%%%%%%%%%%%%

%%%%%%%%%%%%%%%%%%%%%%%%%%%%%%%%%%%%%%%%%%%%%%%%%%%%%%%%%%%%
\subsection{La $q$-conjecture des d\'enominateurs}
D'apr\`es le Lemme~3.5, le
d\'enominateur commun \`a tous les coefficients $\hat{P}_{j,n}(q)$
($j\in\{0,3,5,\dots,A-1\}$) dans l'expression
(\ref{formelineaire}) est de la forme
$D_n(q)=(A-1)!\,q^{\lfloor\alpha n^2+\beta
n+\gamma\rfloor}d_n(1/q)^{A}$, avec $\alpha=-A/8-r^2/2$. Comme
nous l'avons vu  dans la partie pr\'ec\'edente, l'am\'elioration 
donn\'ee par le Th\'eor\`eme~1.1 correspond au gain d'une puissance
de $d_n(1/q)$ dans ce d\'enominateur commun. Or il apparait
clairement dans (\ref{denominatpjchap}) que le d\'enominateur
commun aux $\hat{P}_{j,n}(q)$ pour $j\in\{3,5,\dots,A-1\}$ est de
la forme $\tilde D_n(q)=(A-1)!\,q^{\lfloor\alpha n^2+\beta
n+\gamma\rfloor}d_n(1/q)^{A-1}$, avec $\alpha=-A/8-r^2/2$. Il suffit donc de prouver que cela reste valable pour $\hat{P}_{0,n}(q)$. Ainsi le Th\'eor\`eme~$1.1$ est une cons\'equence du r\'esultat suivant, qui est l'expression de notre $q$-conjecture des d\'enominateurs~:
\begin{theo}
Soient $A$ entier pair et $r\in\N^*$ tel que $A-2r>0$. Soit $\alpha=-A/8-r^2/2$; il existe $\beta$ et $\gamma$ r\'eels ne
d\'ependant que de $A$ et $r$ tels que :
$$q^{\lfloor\alpha n^2+\beta
n+\gamma\rfloor}d_n(1/q)^{A-1}\hat{P}_{0,n}(q)\in\Z\left[\frac{1}{q}\right].$$
\end{theo}

\begin{rem}
\emph{Il est possible, mais inutile dans le cadre de la preuve du
Th\'eor\`eme~1.1 \`a laquelle se limite cet article, de prouver
avec des outils similaires \`a ceux qui suivent qu'en fait~:
\begin{equation*}
(A-1)!\,q^{\lfloor\alpha n^2+\beta
n+\gamma\rfloor}d_n(1/q)^{A-j-1}\hat{P}_{j,n}(q)\in\Z\left[\frac{1}{q}\right]\;\;\forall
j\in\{0,3,5,\dots,A-1\}.
\end{equation*}}
\end{rem}

%%%%%%%%%%%%%%%%%%%%%%%%%%%%%%%%%%%%%%%%%%%%%%%%%%%%%%%%%%%%
\subsection{Une condition suffisante}

En vertu du fait que $\lim_{q\to+\infty}d_n(1/q)=d_n(0)=\pm1$ et
de la formule (\ref{denominatp0chap}) du Lemme~3.5, il suffit pour
prouver le Th\'eor\`eme~4.1 de montrer que
$$d_n(1/q)^{A-1}\hat{P}_{0,n}(q)\in\Z\left[q;\frac{1}{q}\right].$$
D'abord, le troisi\`eme terme de \eqref{P0chap} servant dans le
Lemme~3.2 \`a d\'efinir $\hat{P}_{0,n}(q)$ v\'erifie gr\^ace \`a
(\ref{100}) et  (\ref{*}) :
$$d_n(1/q)^{A-1}\left[\frac{d}{dz}\tilde{P}_{1,n}(z,q)\right]_{z=1}\in\Z\left[q;\frac{1}{q}\right].$$
Il suffit donc de d\'emontrer que
$$d_n(1/q)^{A-1}\left(\tilde{P}_{0,n}(1,q)-q^{-n(r-1)}\tilde{P}_{0,n}(1,1/q)\right)\in\Z\left[q;\frac{1}{q}\right].$$

\medskip
\noindent 
Reprenons donc la d\'efinition (\ref{P0tilde}) de
$\tilde{P}_{0,n}(z,q)$, qui \`a l'aide d'un calcul simple et de
(\ref{relationcj}) permet d'\'ecrire :
\begin{equation}\label{110}
\tilde{P}_{0,n}(1,q)-q^{-n(r-1)}\tilde{P}_{0,n}(1,1/q)=\sum_{s=1}^A\sum_{k=1}^n\frac{1}{(1-q^{-k})^s}V_k,
\end{equation}
o\`u
\begin{equation}\label{Vk}
V_k=\frac{1}{(A-s)!}\left[\frac{d^{A-s}}{du^{A-s}}\sum_{j=k}^nq^{-k(s-1)}e_j(u)-(-1)^sq^{-k}e_{n-j}(u)\right]_{u=1},
\end{equation}
avec comme notation
\begin{eqnarray}
e_j(u)&:=&q^{-j}\tilde{R}_n(q^{-j}u,q)(u-1)^A\nonumber\\
&=&(-1)^{rn}u^{A/2-2+(n-2j)A/2}q^{-r^2n^2/2-nr/2-j(n-j)A/2+j}\nonumber\\
&&\!\times\left(\frac{(q)_{rn}}{(q)_n^r}\right)^2\left(\frac{(q)_n}{(qu^{-1})_j(qu)_{n-j}}\right)^A\frac{(q^{j+1}u^{-1},q^{n+1-j}u)_{rn}}{(q,q)_{rn}}\,\cdot\label{ej}
\end{eqnarray}
La formule \eqref{110} est en fait obtenue en inversant les sommes en $j$ et en $k$ d\'efinies en (\ref{P0tilde}). Cette \'etape d\'eterminante nous permettra d'extraire un facteur de $V_k$ (voir \eqref{alpha} ci-dessous).
Notons que l'on a
\begin{equation}\label{lien-c-e}
\tilde{c}_{s,j,n}(q)=\frac{1}{(A-s)!}\left[q^{-j(s-1)}\frac{d^{A-s}}{du^{A-s}}e_j(u)\right]_{u=1},
\end{equation}
et que $e_j(u)$ est en fait d\'efini pour $A\in2\N$ et $r\in\N$ quelconque.

\medskip

Du Lemme~3.3 et de la relation (\ref{lien-c-e}), on d\'eduit que
pour tout $s\in\{1,\dots,A\}$ et $j\in\{1,\dots,n\}$,
$d_n(1/q)^{A-s}\left[d^{A-s}e_j(u)/du^{A-s}\right]_{u=1}\in\Z\left[q;1/q\right]$,
et par suite que pour tout $k\in\{1,\dots,n\}$,
$d_n(1/q)^{A-s}V_k\in\Z\left[q;1/q\right]$. Ceci ne prouve pas le
Th\'eor\`eme~4.1, mais par d\'efinition de $d_n(q)$, il suffit via
l'\'equation (\ref{110}) de montrer que pour tous
$s\in\{1,\dots,A\}$ et $k\in\{1,\dots,n\}$, on a en fait :
\begin{equation}\label{alpha}
\frac{1}{1-q^{-k}}\,d_n(1/q)^{A-s}V_k\in\Z\left[q;\frac{1}{q}\right].
\end{equation}
C'est ceci qui se r\'ev\`ele difficile, tous les termes
d\'efinissant $V_k$ \'etant importants. En effet, les calculs
effectu\'es avec le logiciel Maple montrent que cela est faux si
l'on remplace $V_k$ d\'efini en (\ref{Vk}) par
$1/(A-s)!\times\left[d^{A-s}e_j(u)/du^{A-s}\right]_{u=1}$ ou
m\^eme $1/(A-s)!\times\sum_{j=k}^nq^{-k(s-1)}
\left[d^{A-s}e_j(u)/du^{A-s}\right]_{u=1}$. 

\medskip

Ecrivons maintenant notre \emph{condition suffisante} nous permettant d'obtenir le Th\'eor\`eme~4.1 et donc le Th\'eor\`eme~1.1~: pour tous $s\in\{1,\dots,A\}$ et $k\in\{1,\dots,n\}$,
\begin{equation}\label{cs}
\frac{1}{1-q^{-k}}\,\frac{d_n(1/q)^{A-s}}{(A-s)!}\left[\frac{d^{A-s}}{du^{A-s}}\sum_{j=k}^ne_j(u)(1-q^{n-2j}u^2)\right]_{u=1}\in\Z\left[q;\frac{1}{q}\right].
\end{equation}
On a en effet le lemme suivant~:
\begin{lem}
Soit $A$ entier pair, et $r\in\N^*$ entier tel que $A-2r>0$. Fixons $s\in\{1,\dots,A\}$ et $k\in\{1,\dots,n\}$. Si
$$
\frac{1}{1-q^{-k}}\,\frac{d_n(1/q)^{A-s}}{(A-s)!}\left[\frac{d^{A-s}}{du^{A-s}}\sum_{j=k}^ne_j(u)(1-q^{n-2j}u^2)\right]_{u=1}\in\Z\left[q;\frac{1}{q}\right],
$$
alors
$$\frac{1}{1-q^{-k}}\,d_n(1/q)^{A-s}V_k\in\Z\left[q;\frac{1}{q}\right].$$
\end{lem}

\medskip

La condition suffisante \eqref{cs} sera d\'emontr\'ee par le lemme cl\'e du paragraphe~4.3.
 
\begin{rem}
\emph{On a ainsi remplac\'e la somme d\'efinissant $V_k$ en
\eqref{Vk} par $\sum_{j=k}^ne_j(u)(1-q^{n-2j}u^2)$. L'int\'er\^et
de cette derni\`ere somme est que, contrairement \`a celle
d\'efinissant $V_k$, elle peut s'\'ecrire comme limite d'une
s\'erie hyperg\'eom\'etrique basique tr\`es bien \'equilibr\'ee,
qui elle-m\^eme pourra s'exprimer (via une transformation due \`a
Andrews donn\'ee au paragraphe~4.4) \`a l'aide d'une somme multiple
$\sum_{\underline{j}}v_{\underline{j}}(u)$, dont chaque terme
poss\`edera la propri\'et\'e souhait\'ee, \`a savoir :
$$\frac{1}{1-q^{-k}}\,\frac{d_n(1/q)^{A-s}}{(A-s)!}\left[\frac{d^{A-s}}{du^{A-s}}v_{\underline{j}}(u)\right]_{u=1}\in\Z\left[q;\frac{1}{q}\right].$$}
\end{rem}

\medskip

\begin{proof}[D\'emonstration du Lemme 4.3]
$V_k$ d\'efini par (\ref{Vk}) peut aussi s'\'ecrire
\begin{multline*}
V_k=\frac{1}{(A-s)!}\left[\frac{d^{A-s}}{du^{A-s}}\sum_{j=k}^n(q^{-k(s-1)}-q^{-k})e_j(u)\right]_{u=1}\\
+\frac{q^{-k}}{(A-s)!}\left[\frac{d^{A-s}}{du^{A-s}}\sum_{j=k}^ne_j(u)-(-1)^se_{n-j}(u)\right]_{u=1}.
\end{multline*}
Or $e_j(u)$ d\'efini par (\ref{ej}) v\'erifie la relation
\begin{equation}\label{relationeju}
e_{n-j}(u)=u^{A-4}q^{n-2j}e_j(1/u).
\end{equation}
Si l'on d\'erive $l$ fois cette \'egalit\'e par rapport \`a $u$,
et l'on sp\'ecialise $u=1$, alors on obtient pour des
$\lambda_i\in\Z$ :
$$\left[\frac{d^{l}}{du^{l}}e_{n-j}(u)\right]_{u=1}=(-1)^lq^{n-2j}\left[\frac{d^{l}}{du^{l}}e_{j}(u)\right]_{u=1}+q^{n-2j}\sum_{i=0}^{l-1}\lambda_i\left[\frac{d^{i}}{du^{i}}e_{j}(u)\right]_{u=1}.$$

\medskip
\noindent
De cela, nous d\'eduisons en posant $l=A-s$ que $V_k$ peut
s'\'ecrire
\begin{equation}\label{140}
V_k=V_{k,1}+V_{k,2}+\frac{q^{-k}}{(A-s)!}\left[\frac{d^{A-s}}{du^{A-s}}\sum_{j=k}^n(1-q^{n-2j})e_{j}(u)\right]_{u=1},
\end{equation}
avec
$$V_{k,1}=(q^{-k(s-1)}-q^{-k})\left[\frac{d^{A-s}}{du^{A-s}}\sum_{j=k}^ne_{j}(u)\right]_{u=1},$$
et pour des $(\alpha_i,\beta_i)\in\Z^2$
$$V_{k,2}=\sum_{i=0}^{A-s-1}q^{\alpha_i}\beta_i\left[\frac{d^{i}}{du^{i}}e_{j}(u)\right]_{u=1}.$$

\medskip
\noindent
Tout d'abord, concernant $V_{k,1}$, pour tout $s\in\{1,\dots,A\}$ on a
$(q^{-k(s-1)}-q^{-k})/(1-q^{-k})\in\mathbb{Z}[1/q]$. Par ailleurs, avec le
Lemme~3.3 et la relation (\ref{lien-c-e}), on a~:
$d_n(1/q)^{A-s}/(A-s)!\times
[d^{A-s}e_j(u)/du^{A-s}]_{u=1}\in\mathbb{Z}[q;1/q]$. On en d\'eduit que~:
\begin{equation}\label{142}
\frac{1}{1-q^{-k}}\,d_n(1/q)^{A-s}V_{k,1}\in\mathbb{Z}\left[q;\frac{1}{q}\right].
\end{equation}
Mais on a aussi gr\^ace au Lemme~3.3 et \`a la relation
(\ref{lien-c-e}) $d_n(1/q)^{A-s-1}V_{k,2}\in\Z[q;1/q]$, donc \`a
fortiori
\begin{equation}\label{146}
\frac{1}{1-q^{-k}}\,d_n(1/q)^{A-s}V_{k,2}\in\mathbb{Z}\left[q;\frac{1}{q}\right].
\end{equation}
Enfin on peut \'ecrire~:
\begin{multline}\label{150}
\frac{1}{(A-s)!}\left[\frac{d^{A-s}}{du^{A-s}}\sum_{j=k}^n(1-q^{n-2j})e_{j}(u)\right]_{u=1}\\
=\frac{1}{(A-s)!}\left[\frac{d^{A-s}}{du^{A-s}}\sum_{j=k}^n(1-q^{n-2j}u^2)e_{j}(u)\right]_{u=1}+V_{k,3},
\end{multline}
avec $V_{k,3}$ exactement du m\^eme type que $V_{k,2}$, donc
v\'erifiant
\begin{equation}\label{148}
\frac{1}{1-q^{-k}}\,d_n(1/q)^{A-s}V_{k,3}\in\mathbb{Z}\left[q;\frac{1}{q}\right].
\end{equation}
Finalement les expressions (\ref{140}) et (\ref{150}) permettent
de conclure, via (\ref{142}), (\ref{146}) et (\ref{148}).
\end{proof}

%%%%%%%%%%%%%%%%%%%%%%%%%%%%%%%%%%%%%%%%%%%%%%%%%%%%%%%%%%%%
\subsection{Le lemme cl\'e}

D'apr\`es le Lemme~4.3, il suffit pour prouver le Th\'eor\`eme~4.1 de d\'emontrer la condition suffisante \eqref{cs}. Celle-ci est elle-m\^eme cons\'equence du lemme plus g\'en\'eral suivant~:
\begin{lem}
Soient $l$, $r$, $A$, $k$ et $n$ des entiers positifs, avec $A$ pair et $k\in\{1,\dots,n\}$. On a~:
$$\frac{1}{1-q^{-k}}\,\frac{d_n\left(1/q\right)^{l}}{l!}\left[\frac{d^{l}}{du^{l}}\sum_{j=k}^n(1-q^{n-2j}u^2)e_j(u)\right]_{u=1}\in\mathbb{Z}\left[q;\frac{1}{q}\right],$$
o\`u $e_j(u)$ (d\'efini en \eqref{ej}) s'\'ecrit~:
\begin{multline*}
e_j(u)=u^{A/2-2+(n-2j)A/2}q^{-r^2n^2/2-nr/2-j(n-j)A/2+j}\\
\times(-1)^{rn}\left(\frac{(q)_{rn}}{(q)_n^r}\right)^2\left(\frac{(q)_n}{(qu^{-1})_j(qu)_{n-j}}\right)^A\frac{(q^{j+1}u^{-1},q^{n+1-j}u)_{rn}}{(q,q)_{rn}}\,\cdot
\end{multline*}
\end{lem}

\noindent Pour d\'emontrer ce lemme, nous avons besoin d'une identit\'e
g\'en\'erale de transformation de s\'eries hyperg\'eom\'etriques
basiques, ce qui fait l'objet du paragraphe $4.4$. La
d\'emonstration du Lemme~4.5 sera donn\'ee dans le paragraphe~$4.5$. 

\medskip

Afin d'illustrer le Lemme~4.5 lorsque $l=0$, \'ecrivons (au signe pr\`es) $\sum_{j=k}^n(1-q^{n-2j})e_j(1)$, dans les cas particuliers $(A,r)=(0,0)$,  $(A,r)=(2,0)$, $(A,r)=(0,1)$ et $(A,r)=(2,1)$ respectivement (en posant $p=1/q$)~:

\begin{eqnarray}
&&\hskip -1.5cm\sum_{j=k}^n(1-q^{n-2j})q^{j}=(1-p^{k})p^{-n}\frac{1-p^{n-k+1}}{1-p}\label{0,0},\\
&&\hskip -1.5cm\sum_{j=k}^n(1-q^{n-2j})q^{j(j-n+1)}\qbi{n}{j}{q}^2
=(1-p^{k})p^{-k(n-k+1)}\qbi{n}{k}{p}\qbi{n-1}{k-1}{p}\label{2,0},\\
&&\hskip -1.5cm\sum_{j=k}^n(1-q^{n-2j})q^{j-n(n+1)/2}\qbi{n+j}{n}{q}\qbi{2n-j}{n}{q}\nonumber\\
&&\hskip 2.5cm=(1-p^{k})p^{-n(n+1)/2}\qbi{n+k}{n}{p}\qbi{2n+1-k}{n+1}{p}\label{0,1},
\end{eqnarray}
\begin{multline}\label{2,1}
\sum_{j=k}^n(1-q^{n-2j})q^{j^2+j-nj-n(n+1)/2}\qbi{n+j}{n}{q}\qbi{2n-j}{n}{q}\qbi{n}{j}{q}^2\\
=(1-p^{k})p^{k^2-kn-k-n(n-1)/2}\sum_{l=0}^{n-k}(-1)^lp^{l(2k+l-1)/2}\\
\times\qbi{n+k}{k+l}{p}\qbi{n-l-1}{k-1}{p}\qbi{2n-k-l}{k,n-k,n-k-l}{p}.
\end{multline}
L'\'egalit\'e \eqref{0,0} est obtenue par un calcul direct, alors que 
\eqref{2,0}-\eqref{2,1} sont obtenues par des sp\'ecialisations de la transformation de Watson finie \cite[Appendix III, (III.18)]{GR}. Ces identit\'es montrent bien que  $\sum_{j=k}^n(1-q^{n-2j})\,e_j(1)\in(1-q^{-k})\,\mathbb{Z}\left[q;1/q\right]$. En particulier, \eqref{2,1} illustre le fait que ceci peut \^etre visible malgr\'e un membre de droite compliqu\'e.

\medskip
 Lorsque l'on divise toutes ces \'egalit\'es par $-(1-q)=(1-p)p^{-1}$ et l'on fait tendre $q$ vers 1, on trouve respectivement les identit\'es hyperg\'eom\'etriques suivantes, dont les trois premi\`eres ont \'et\'e donn\'ees dans \cite[p. 62]{KR}~: 
\begin{eqnarray*}
&&-\sum_{j=k}^n(n-2j)=k(n-k+1),\\
&&-\sum_{j=k}^n(n-2j)\bi{n}{j}^2
=k\bi{n}{k}\bi{n-1}{k-1},\\
&&-\sum_{j=k}^n(n-2j)\bi{n+j}{n}\bi{2n-j}{n}
=k\bi{n+k}{k}\bi{2n+1-k}{n+1},\\
&&-\sum_{j=k}^n(n-2j)\bi{n+j}{n}\bi{2n-j}{n}\bi{n}{j}^2\\
&&\hskip 2cm=k\sum_{l=0}^{n-k}(-1)^l\bi{n+k}{k+l}\bi{n-l-1}{k-1}\bi{2n-k-l}{k,n-k,n-k-l}.
\end{eqnarray*}

%%%%%%%%%%%%%%%%%%%%%%%%%%%%%%%%%%%%%%%%%%%%%%%%%%%%%%%%%%%%
\subsection{Une transformation g\'en\'erale de $q$-s\'eries}

Nous avons maintenant besoin, pour montrer le Lemme~4.5, d'une identit\'e de transformation due \`a Andrews \cite{An75, An84}
entre une somme simple tr\`es bien \'equilibr\'ee et une somme
multiple, qui est le $q$-analogue du Th\'eor\`eme 8 de \cite{KR}.
Cette identit\'e a \'et\'e prouv\'ee d'abord dans \cite{An75},
avant d'\^etre vue dans \cite{An84} comme une cons\'equence
directe du lemme de Bailey \cite{Ba}, qui est lui-m\^eme un outil
\'el\'ementaire et tr\`es efficace pour prouver des identit\'es de
$q$-s\'eries.
\begin{theo}[Andrews]
Pour tous entiers $m\geq 0$ et $N\geq 0$, pour tous complexes $a$,
$b_1,c_1,\dots,b_{m+1},c_{m+1}$, on a :
\begin{multline}\label{andrews}
\sum_{k=0}^N\frac{1-aq^{2k}}{1-a}\frac{(a,b_1,c_1,\dots,b_{m+1},c_{m+1},q^{-N})_k}
{(q,aq/b_1,aq/c_1,\dots,aq/b_{m+1},aq/c_{m+1},aq^{N+1})_k}
\\
\times\left(\frac{a^{m+1}q^{m+1+N}}{b_1c_1\dots b_{m+1}c_{m+1}}\right)^k\\
=\frac{(aq,aq/b_{m+1}c_{m+1})_N}{(aq/b_{m+1},aq/c_{m+1})_N}\sum_{0\leq
l_1\leq\dots\leq l_{m}\leq N}
\frac{a^{l_1+\dots+l_{m-1}}q^{l_1+\dots+l_{m}}}{(b_2c_2)^{l_1}\dots(b_mc_m)^{l_{m-1}}}\\
\times\frac{(q^{-N})_{l_m}}{(b_{m+1}c_{m+1}q^{-N}/a)_{l_m}}\prod_{i=1}^m\frac{(b_{i+1},c_{i+1})_{l_i}}{(aq/b_i,aq/c_i)_{l_i}}\frac{(aq/b_ic_i)_{l_i-l_{i-1}}}{(q)_{l_i-l_{i-1}}}\cdot
\end{multline}
\end{theo}
Cette \'egalit\'e justifie notre intention aux paragraphes pr\'ec\'edants de nous
ramener \`a une s\'erie tr\`es bien \'equilibr\'ee. En effet, lorsque pour tout $i\in\{1,\dots,m+1\}$, $aq/b_i\neq q^{-N}$ et $aq/c_i\neq q^{-N}$, le membre de gauche de \eqref{andrews} peut s'\'ecrire en terme de
s\'erie hyperg\'eom\'etrique basique tr\`es bien \'equilibr\'ee~:
\begin{equation*}\label{hyperg}
{}_{2m+6}\phi_{2m+5}\!\left[\begin{matrix}a,q\sqrt{a},-q\sqrt{a},b_1,c_1,\dots,b_{m+1},c_{m+1},q^{-N}\\
\sqrt{a},-\sqrt{a},aq/b_1,aq/c_1,\dots,aq/b_{m+1},aq/c_{m+1},aq^{N+1}\end{matrix};q,z\right],
\end{equation*}
avec $z=a^{m+1}q^{m+1+N}/b_1c_1\dots b_{m+1}c_{m+1}$. Cependant, si l'un des facteurs $b_i$ ou $c_i$, par exemple $b_1$, s'\'ecrit $aq^{1+N}$, alors la s\'erie hyperg\'eom\'etrique ci-dessus n'est plus une somme finie. En effet, le facteur $(q^{-N})_k$ du num\'erateur, garantissant cette finitude, se simplifie avec le d\'enominateur $(aq/b_1)_k$. Le membre de gauche de \eqref{andrews}, qui est toujours une somme finie, peut alors \^etre vu comme limite d'une s\'erie hyperg\'eom\'etrique basique tr\`es bien \'equilibr\'ee~:
\begin{equation*}
\lim_{\delta\to 0}\;{}_{2m+6}\phi_{2m+5}\!\left[\begin{matrix}a,q\sqrt{a},-q\sqrt{a},aq^{1+N+\delta},c_1,\dots,q^{-N}\\
\sqrt{a},-\sqrt{a},q^{-N-\delta},aq/c_1,\dots,aq^{N+1}\end{matrix};q,zq^{-\delta}\right].
\end{equation*}

\medskip
\begin{rem}
\emph{Dans le cas $m=1$, l'\'egalit\'e \eqref{andrews} est exactement la transformation finie de Watson d\'ej\`a mentionn\'ee dans le paragraphe pr\'ec\'edent et utile pour prouver les identit\'es \eqref{2,0}-\eqref{2,1}. Ceci indique bien que \eqref{andrews} pourrait \^etre l'ingr\'edient manquant pour d\'emontrer le Lemme~4.5.}
\end{rem}

%%%%%%%%%%%%%%%%%%%%%%%%%%%%%%%%%%%%%%%%%%%%%%%%%%%%%%%%%%%%
\subsection{D\'emonstration du Lemme~4.5}

Posons
$v_j(u):=e_j(u)(1-q^{n-2j}u^2)$, de sorte que le crochet de
l'expression que l'on \'etudie dans le Lemme~4.5 peut s'\'ecrire
$$
\sum_{j=k}^n(1-q^{n-2j}u^2)e_j(u)=\sum_{j=k}^nv_j(u)=\sum_{j=0}^{n-k}v_{k+j}(u).$$
Ceci se transforme en
\begin{multline}\label{exprepourandrews'}
v_k(u)\times\sum_{j=0}^{n-k}p^{j((A/2-r)n+A/2-1)}\times\frac{1-ap^{2j}}{1-a}\\
\times\frac{(a,p,p^{rn+k+1}u,p^{k-n}u,\dots,p^{k-n}u,p^{k+1}u^2,p^{k-n};p)_j}{(p,a,p^{k-rn-n}u,p^{k+1}u,\dots,p^{k+1}u,p^{k-n},p^{k+1}u^2;p)_j},
\end{multline}
 avec $p=1/q$ et $a=p^{2k-n}u^2$, les facteurs $p^{k-n}u$ au num\'erateur et $p^{k+1}u$ au d\'enominateur apparaissant $A+1$ fois. On peut maintenant appliquer la formule d'Andrews
(\ref{andrews}) \`a (\ref{exprepourandrews'}) en rempla\c cant $q$
par $p$, et en posant $N=n-k$, $m=A/2+1$, $a=p^{2k-n}u^2$,
$b_1=p^{rn+k+1}u$,
$c_1=b_2=c_2=\dots=b_{A/2}=c_{A/2}=b_{A/2+1}=c_{A/2+2}=p^{k-n}u$,
$c_{A/2+1}=p^{k+1}u^2$ et $b_{A/2+2}=p$. Ceci donne
comme expression alternative pour (\ref{exprepourandrews'}) :
\begin{multline*}
v_k(u)\times\frac{(p^{2k+1-n}u^2,p^{k}u;p)_{n-k}}{(p^{2k-n}u^2,p^{k+1}u;p)_{n-k}}\\
\times\sum_{0\leq l_1\leq\dots\leq l_{A/2+1}\leq
n-k}\frac{(p^{2k-n}u^2)^{l_1+\dots+l_{A/2}}p^{l_1+\dots+l_{A/2+1}}}{(p^{2k-2n}u^2)^{l_1+\dots+l_{A/2-1}}(p^{2k-n+1}u^3)^{l_{A/2}}}\\
\times\frac{(p^{-rn},p^{k-n}u,p^{k-n}u;p)_{l_1}}{(p,p^{k-n-rn}u,p^{k+1}u;p)_{l_1}}\times\prod_{i=1}^{A/2}\frac{(p^{k-n}u,p^{k-n}u;p)_{l_i}(p^{n+1};p)_{l_i-l_{i-1}}}{(p^{k+1}u,p^{k+1}u;p)_{l_i}(p;p)_{l_i-l_{i-1}}}\\
\times\frac{(p^{k+1}u^2;p)_{l_{A/2}}}{(p^{k-n}u;p)_{l_{A/2}}}\times\frac{(p^{k-n},p,p^{k-n}u;p)_{l_{A/2+1}}}{(p^{1-n}/u,p^{k+1}u,p^{k-n};p)_{l_{A/2+1}}}\times\frac{(1/u;p)_{l_{A/2+1}-l_{A/2}}}{(p;p)_{l_{A/2+1}-l_{A/2}}}\cdot
\end{multline*}

\medskip
\noindent
En posant $\underline{l}:=(l_1,\dots,l_{A/2+1})$, on obtient
apr\`es maintes transformations \'el\'ementaires de $p$-factoriels
($p=1/q$)~:
\begin{equation}\label{expreapresandrewsbriques}
\sum_{j=k}^n(1-q^{n-2j}u^2)e_j(u)=(1-p^ku)\sum_{0\leq l_1\leq\dots\leq l_{A/2+1}\leq n-k}w_{\underline{l}}(u),
\end{equation}
o\`u pour $\xi(\underline{l})$, $\varphi(\underline{l})$ et
$\psi(\underline{l})$ \'el\'ements de $\Z$,
\begin{multline}\label{expreenbriques}
w_{\underline{l}}(u):=
(-1)^{\xi(\underline{l})}p^{\varphi(\underline{l})}u^{\psi(\underline{l})}\times
R_1(n,k,n,p,u)\times\prod_{i=2}^{A/2}\qbi{n+l_{i}-l_{i-1}}{k+l_i}{p}\\
\times\qbi{l_{A/2+1}}{l_1,l_2-l_1,\dots,l_{A/2+1}-l_{A/2}}{p}\times\prod_{i=1}^{A/2}R_0(0,n+1,k+l_i,p,1/u)\\
\times\prod_{i=2}^{A/2}R_0(0,n-k-l_{i-1}+1,0,p,1/u)\times R_0(0,k+l_{i}+1,0,p,u)\\
\times R_1(n,n-l_1-k,n-l_1,p,1/u)\times R_2(n,k,l_{A/2+1},l_{A/2},p,u),
\end{multline}
et pour $\alpha\leq\gamma<\beta$ entiers et $|q|\neq 1$
\begin{equation}\label{briqueRtilde}
R_0(\alpha,\beta,\gamma,q,u):=\frac{(1-u)\,(q)_{\beta-\alpha-1}}{(q^{\alpha-\gamma}u)_{\beta-\alpha}}=\frac{(q)_{\beta-\alpha-1}}{(q^{\alpha-\gamma}u)_{\gamma-\alpha}(qu)_{\beta-\gamma-1}},
\end{equation}
puis pour $n\geq j\geq i\geq0$ entiers et $|q|\neq 1$
\begin{equation}\label{briqueR1}
R_1(n,i,j,q,u):=\frac{(q)_{rn}}{(q)_n^r}\times\frac{(q^{n+i-j+1}u)_{rn-n+j}}{(q)_{rn-n+j}},
\end{equation}
et enfin pour $0\leq m_2\leq m_1\leq n-k$ entiers et $|q|\neq 1$
\begin{multline}\label{briqueR2}
R_2(n,k,m_1,m_2,q,u):=\frac{(q)_n}{(qu)_n}\frac{(1/u)_{m_1-m_2}}{(qu^2)_{k-1}}\frac{(qu^2)_{k+m_2}}{(qu)_{k+m_1}}\frac{(qu)_{n-m_1-1}}{(q/u)_{n-k-m_1}}\cdot
\end{multline}

\medskip
\begin{rem}
\emph{Les fractions $R_0/(1-u)$, $R_1$ et $R_2$ d\'efinies \`a l'aide de \eqref{briqueRtilde}-\eqref{briqueR2} sont des
$q$-analogues de certaines briques \'el\'ementaires et sp\'eciales de
\cite{KR}.} 
\end{rem}

\medskip
L'expression (\ref{expreapresandrewsbriques}) indique que pour
d\'emontrer le Lemme~4.5, il nous suffit maintenant, en vertu de
la formule de Leibniz de d\'erivation d'un produit, de prouver que pour tout entier $l$~:
$$\frac{d_n\left(p\right)^{l}}{l!}\left[\frac{d^{l}}{du^{l}}w_{\underline{l}}(u)\right]_{u=1}\in\mathbb{Z}\left[p;\frac{1}{p}\right]=\mathbb{Z}\left[q;\frac{1}{q}\right].$$
Mais gr\^ace \`a la formule de Leibniz de nouveau et \`a
l'expression (\ref{expreenbriques}), ceci est une cons\'equence du
lemme suivant~:
\begin{lem}
Soit $|q|\neq 1$. Pour tous entiers $l\geq0$,
$\alpha\leq\gamma<\beta$, $n\geq j\geq
i\geq0$ et $0\leq m_2\leq m_1\leq n-k$, on a :
\begin{eqnarray}
&&\frac{d_{\beta-\alpha-1}(q)^l}{l!}\left[\frac{d^l}{du^l}R_0(\alpha,\beta,\gamma,q,u)\right]_{u=1}\in\Z\left[q,\frac{1}{q}\right],\label{Rtilde}\\
&&\frac{d_{n}(q)^l}{l!}\left[\frac{d^l}{du^l}R_1(n,i,j,q,u)\right]_{u=1}\in\Z\left[q,\frac{1}{q}\right],\label{R1}\\
&&\frac{d_{n}(q)^l}{l!}\left[\frac{d^l}{du^l}R_2(n,k,m_1,m_2,q,u)\right]_{u=1}\in\Z\left[q,\frac{1}{q}\right]\label{R2}.
\end{eqnarray}
De plus, ces propri\'et\'es restent valables lorsque $u$ est
remplac\'e par $u^e$, $e\in\Z$ (en particulier par $1/u$).
\end{lem}
\begin{proof}[D\'emonstration]
Tout d'abord, le fait que les propri\'et\'es
(\ref{Rtilde})-(\ref{R2}) restent valables lorsque $u$ est
remplac\'e par $u^e$, $e\in\Z$, est une simple cons\'equence de la
d\'erivation des fonctions
compos\'ees.

\medskip

\noindent Commen\c cons donc par d\'emontrer (\ref{Rtilde}), en
d\'ecomposant en \'el\'ements simples la fraction rationnelle (en
$u$)~:
$$
R_0(\alpha,\beta,\gamma,q,u)=\sum_{j=\alpha\atop
j\neq\gamma}^{\beta-1}a_j\frac{1-q^{j-\gamma}}{1-q^{j-\gamma}u},
$$
avec $$\displaystyle
a_j=\frac{(1-q)\dots(1-q^{\beta-\alpha-1})}{(1-q^{\alpha-j})\dots(1-q^{-1})(1-q)\dots(1-q^{\beta-j-1})},$$
qui se r\'ecrit
$$a_j=(-1)^{j-\alpha} q^{(j-\alpha)(j-\alpha+1)/2}\qbi{\beta-\alpha-1}{j-\alpha}{q}\in\Z\left[q\right].$$
En d\'erivant, on obtient donc 
$$
\frac{1}{l!}\frac{d^l}{du^l}R_0(\alpha,\beta,\gamma,q,u)=\sum_{j=\alpha\atop
j\neq\gamma}^{\beta-1}q^{(j-\gamma)l}a_j\frac{1-q^{j-\gamma}}{(1-q^{j-\gamma}u)^{l+1}},
$$
ce qui en prenant $u=1$ prouve (\ref{Rtilde}).\\

\medskip

\noindent D\'emontrons ensuite (\ref{R1}). Comme $R_1(n,i,j,q,u)$ est un polyn\^ome en $u$ de degr\'e $(r-1)n+j$, on a $\frac{d^l}{du^l}R_1(n,i,j,q,u)=0$ si $l>(r-1)n+j$. Supposons donc que $l\leq(r-1)n+j$. On a alors~:
\begin{multline*}
\frac{1}{l!}\frac{d^l}{du^l}R_1(n,i,j,q,u)=\sum_{n+i-j+1\leq
f_1<\dots< f_l\leq
rn+i}(-1)^lq^{f_1+\dots+f_l}\\
\times\frac{R_1(n,i,j,q,u)}{(1-q^{f_1}u)\dots(1-q^{f_l}u)}\cdot
\end{multline*}
En posant $g_m:=f_m-(n+i-j)$ pour tout $m\in\{1,\dots,l\}$, ceci peut s'\'ecrire lorsque $u=1$~:
$$
\frac{1}{l!}\left[\frac{d^l}{du^l}R_1(n,i,j,q,u)\right]_{u=1}=\sum_{1\leq
g_1<\dots< g_l\leq rn-n+j}(-1)^lq^{(n+i-j)l+g_1+\dots+g_l}\,\tilde{R}_1(q),
$$
avec 
\begin{equation}\label{R1tildeq}
\tilde{R}_1(q):=\frac{(q)_{rn}}{(q)_n^r}
\frac{(q^{n+i-j+1})_{rn-n+j}}{(q)_{rn-n+j}}\frac{1}{(1-q^{n+i-j+g_1})\dots(1-q^{n+i-j+g_l})}\cdot
\end{equation}
Pour prouver (\ref{R1}), il suffit donc de d\'emontrer que
\begin{equation}\label{R1tilde}
d_{n}(q)^l\tilde{R}_1(q)\in\Z\left[q\right].
\end{equation}
Pour cela, faisons quelques rappels sur les polyn\^omes
cyclotomiques. Pour $t\in\N$, le $t$-i\`eme polyn\^ome
 cyclotomique est $\phi_t(x):=\prod_{k\wedge t=1,k\leq t}(x-\mbox{e}^{2ik\pi/t})$,
avec comme premi\`ere propri\'et\'e
 le fait que $\phi_t(x)\in\Z[x]$. Ensuite, on peut montrer que
\begin{equation}\label{1000}
q^n-1=\prod_{d|n}\phi_d(q),
\end{equation}
 ce qui permet de d\'eduire d'une part que
\begin{equation}\label{dnaveccyclot}
d_n(q)=\prod_{t=1}^n\phi_t(q),
\end{equation}
et d'autre part l'ordre de divisibilit\'e suivant dans $\Z[q]$ :
\begin{equation}\label{ordreqfacto}
ord_{\phi_t(q)}\left(\prod_{l=1}^n\frac{q^l-1}{q-1}\right)=\left\{\begin{array}{l}0\;\;\mbox{si}\;\;t=1\\
\left\lfloor \frac{n}{t}\right\rfloor\;\;\mbox{si}\;\;t\geq
2.\end{array}\right.
\end{equation}
Maintenant, comme $\tilde{R}_1(q)$ est une fraction rationnelle en $q$, dont
num\'erateur et d\'enominateur s'expriment au signe pr\`es en
produits de fonctions du type $\phi_t(q)$ (de par (\ref{1000})),
il suffit pour prouver (\ref{R1tilde}) de voir que
\begin{equation}\label{R1tildeaveccyclot}
\forall\;t\in\N,\;\;ord_{\phi_t(q)}\left(d_{n}(q)^l\tilde{R}_1(q)\right)\geq0.
\end{equation}
Nous distinguons pour cela trois cas.
\begin{itemize}
\item  Si $t=1$, alors on a en vertu de (\ref{dnaveccyclot}) et
(\ref{ordreqfacto})
$$ord_{\phi_1(q)}\left(d_{n}(q)^l\tilde{R}_1(q)\right)=l+ord_{\phi_1(q)}\left(\tilde{R}_1(q)\right)=l-l=0.$$
\item Si $2\leq t\leq n$, alors (\ref{dnaveccyclot}) et
(\ref{ordreqfacto}) donnent
\begin{multline*}
ord_{\phi_t(q)}\left(d_{n}(q)^l\tilde{R}_1(q)\right)=\sum_{m=1}^l\left(1-ord_{\phi_t(q)}(1-q^{g_m+n+i-j})\right)\\
+\left(\left\lfloor \frac{rn}{t}\right\rfloor-r\left\lfloor
\frac{n}{t}\right\rfloor\right)+\left(\left\lfloor
\frac{rn+i}{t}\right\rfloor-\left\lfloor
\frac{n+i-j}{t}\right\rfloor-\left\lfloor
\frac{rn-n+j}{t}\right\rfloor\right),
\end{multline*}
o\`u il est facile de voir que le terme \`a l'int\'erieur de
chaque parenth\`ese est positif ou nul, ce qui permet de conclure.
\item Si $t>n$, on a dans un premier temps par
(\ref{dnaveccyclot})
$$ord_{\phi_t(q)}\left(d_{n}(q)^l\tilde{R}_1(q)\right)=ord_{\phi_t(q)}\left(\tilde{R}_1(q)\right).$$
D'autre part, puisque $1\leq g_1<\dots< g_l\leq rn-n+j$, on peut \'ecrire en utilisant (\ref{1000})
\begin{multline*}
\sum_{m=1}^lord_{\phi_t(q)}\left(1-q^{g_m+n+i-j}\right)\leq\sum_{g=n+i-j+1}^{rn+i}ord_{\phi_t(q)}\left(1-q^{g}\right)\\
=\sum_{g=n+i-j+1\atop
t|g}^{rn+i}ord_{\phi_t(q)}\left(1-q^{g}\right)=\sum_{g=n+i-j+1\atop
t|g}^{rn+i}1=\left\lfloor \frac{rn+i}{t}\right\rfloor-\left\lfloor
\frac{n+i-j}{t}\right\rfloor\\
=ord_{\phi_t(q)}\left((q^{n+i-j+1})_{rn-n+j}\right).
\end{multline*}
Donc finalement, puisque $n\geq j$ et $t>n$, on a en reprenant \eqref{R1tildeq}
\begin{multline*}
ord_{\phi_t(q)}\left(d_{n}(q)^l\tilde{R}_1(q)\right)\geq
\left\lfloor \frac{rn}{t}\right\rfloor-r\left\lfloor \frac{n}{t}\right\rfloor-\left\lfloor
\frac{rn-n+j}{t}\right\rfloor\\
=\left\lfloor \frac{rn}{t}\right\rfloor-\left\lfloor
\frac{rn-n+j}{t}\right\rfloor\geq 0.
\end{multline*}
\end{itemize}

\medskip

\noindent Il ne reste plus qu'\`a prouver (\ref{R2}). Supposons dans un premier temps que $m_1>m_2$. On
peut alors \'ecrire en reprenant la d\'efinition \eqref{briqueR2}
$$R_2(n,k,m_1,m_2,q,u)=(1-u^{-1})R_3(q,u),$$
o\`u 
$$
R_3(q,u):=\frac{(q)_n}{(qu)_n}\frac{(q/u)_{m_1-m_2-1}}{(qu^2)_{k-1}}\frac{(qu^2)_{k+m_2}}{(qu)_{k+m_1}}\frac{(qu)_{n-m_1-1}}{(q/u)_{n-k-m_1}}\cdot
$$
Ceci montre, \`a l'aide de la d\'erivation des produits par
Leibniz, qu'il suffit pour prouver (\ref{R2}) de montrer que
\begin{equation}\label{R3}
\frac{d_{n}(q)^{l+1}}{l!}\left[\frac{d^l}{du^l}R_3(q,u)\right]_{u=1}\in\Z\left[q,\frac{1}{q}\right].
\end{equation}
Mais, sachant que $m_1-m_2-1\leq n$ et $k+m_2\leq k+m_1\leq n$, la
propri\'et\'e de Leibniz et la d\'erivation des fonctions
compos\'ees permettent de voir qu'il suffit pour montrer
(\ref{R3}) de prouver que pour tous $1\leq h_1\leq\dots\leq
h_l\leq n$
\begin{equation}\label{R4}
R_4(q):=d_{n}(q)^{l+1}\frac{R_3(q,1)}{(1-q^{h_1})\dots(1-q^{h_l})}\in\Z\left[q\right].
\end{equation}
On distingue les trois m\^emes cas que pr\'ec\'edemment. Si $t=1$ ou $t>n$, on voit par (\ref{dnaveccyclot}) et
(\ref{ordreqfacto}) que $ord_{\phi_t(q)}\left(R_4(q)\right)=0$. Dans le cas o\`u $2\leq t\leq n$, on \'ecrit
\begin{multline*}
ord_{\phi_t(q)}\left(R_4(q)\right)=1+\sum_{m=1}^l\left(1-ord_{\phi_t(q)}(1-q^{h_m})\right)+ord_{\phi_t(q)}\left(R_3(q,1)\right)\\
\geq \;1+ord_{\phi_t(q)}\left(R_3(q,1)\right)
=1+\left\lfloor
\frac{m_1-m_2-1}{t}\right\rfloor-\left\lfloor \frac{k-1}{t}\right\rfloor\\
+\left\lfloor\frac{k+m_2}{t}\right\rfloor-\left\lfloor \frac{k+m_1}{t}\right\rfloor+\left\lfloor\frac{n-m_1-1}{t}\right\rfloor-\left\lfloor
\frac{n-k-m_1}{t}\right\rfloor\\
=\;
1+\left(\left\lfloor\frac{n-m_1-1}{t}\right\rfloor-\left\lfloor
\frac{n-k-m_1}{t}\right\rfloor-\left\lfloor \frac{k-1}{t}\right\rfloor\right)\\
-\left(\left\lfloor \frac{k+m_1}{t}\right\rfloor-\left\lfloor
\frac{k+m_2}{t}\right\rfloor-\left\lfloor \frac{m_1-m_2}{t}\right\rfloor\right)\\
-\left(\left\lfloor \frac{m_1-m_2}{t}\right\rfloor-\left\lfloor
\frac{m_1-m_2-1}{t}\right\rfloor\right).
\end{multline*}
Or l'identit\'e \'el\'ementaire $0\leq\left\lfloor
a+b\right\rfloor-\left\lfloor a\right\rfloor-\left\lfloor b\right
\rfloor\leq 1$ permet d'affirmer que la premi\`ere parenth\`ese de
cette expression est sup\'erieure ou \'egale \`a 0, et les deux
suivantes sont entre 0 et 1. Ceci montre que si
$ord_{\phi_t(q)}\left(R_4(q)\right)<0$, alors les deux derni\`eres
parenth\`eses doivent valoir 1, ce qui implique que $t$ divise
$m_1-m_2$. Mais en \'ecrivant ces deux derni\`eres parenth\`eses
de deux autres fa\c cons, on voit que
$ord_{\phi_t(q)}\left(R_4(q)\right)<0$ impliquerait aussi le fait
que $t$ divise $k+m_2+1$ et $k+m_1$; donc $t$ diviserait
$1=(m_1-m_2)+(k+m_2+1)-(k+m_1)$,
ce qui est absurde car $t\geq2$. Donc $ord_{\phi_t(q)}\left(R_4(q)\right)\geq 0$.\\

\medskip

\noindent D\'emontrons enfin (\ref{R2}) dans le cas
$m_1=m_2$.  De m\^eme que ci-dessus, il suffit de montrer que pour tous $1\leq h_1\leq\dots\leq h_l\leq n$
$$R_5(q):=d_{n}(q)^{l}\frac{R_2(n,k,m_1,m_1,q,1)}{(1-q^{h_1})\dots(1-q^{h_l})}\in\Z[q].$$
Ici aussi si $t=1$ ou $t>n$, alors par (\ref{dnaveccyclot}) et
(\ref{ordreqfacto}), $ord_{\phi_t(q)}\left(R_5(q)\right)=0$. Pour $2\leq t\leq n$, l'ordre de multiplicit\'e $ord_{\phi_t(q)}\left(R_5(q)\right)$ vaut
$$\begin{aligned}
&\sum_{m=1}^l\left(1-ord_{\phi_t(q)}(1-q^{h_m})\right)+ord_{\phi_t(q)}\left(R_2(n,k,m_1,m_1,q,1)\right)\\
&\hskip 1cm\quad\geq\; ord_{\phi_t(q)}\left(R_2(n,k,m_1,m_1,q,1)\right)\\
&\hskip 1cm\quad=\;\left\lfloor\frac{n-m_1-1}{t}\right\rfloor-\left\lfloor\frac{k-1}{t}\right\rfloor-\left\lfloor
\frac{n-k-m_1}{t}\right\rfloor\,\geq\,0\,.
\end{aligned}$$

\end{proof}

%%%%%%%%%%%%%%%%%%%%%%%%%%%%%%%%%%%%%%%%%%%%%%%%%%%%%%%%%%%%%%%%%%%%%%%%%%%%

\vspace{1cm}

\noindent Fr\'ed\'eric Jouhet,\\
Universit\'e de Lyon, Universit\'e Lyon I, Institut Camille Jordan, UMR 5208,\\
43, bd du 11 Novembre 1918, 69622 Villeurbanne Cedex, France \\
\texttt{jouhet@math.univ-lyon1.fr}

\vspace{0.5cm}

\noindent Elie Mosaki,\\
Universit\'e de Lyon, Universit\'e Lyon I, Institut Camille Jordan, UMR 5208,\\
43, bd du 11 Novembre 1918, 69622 Villeurbanne Cedex, France \\
\texttt{mosaki@math.univ-lyon1.fr}

\end{document}